%
%

\input amstex

\magnification 1200
\loadmsbm
\parindent 0 cm

\define\nl{\bigskip\item{}}
\define\snl{\smallskip\item{}}
\define\inspr #1{\parindent=20pt\bigskip\bf\item{#1}}
\define\iinspr #1{\parindent=27pt\bigskip\bf\item{#1}}
\define\einspr{\parindent=0cm\bigskip}

\define\ot{\otimes}

\define\tr{\triangleright}
\define\tl{\triangleleft}


\input amssym
\input amssym.def

\centerline{\bf Bicrossproducts of multiplier Hopf algebras}
\bigskip\bigskip
\centerline{\it L.\ Delvaux$^{\text{1}}$, A.\ Van Daele$^{\text{2}}$ and S.H.\ Wang$^{\text{3}}$}
\bigskip\bigskip\bigskip
{\bf Abstract} 
\bigskip 
In this paper, we generalize Majid's bicrossproduct construction. We start with a pair $(A,B)$ of two regular multiplier Hopf algebras. We assume that $B$ is a right $A$-module algebra and that $A$ is a left $B$-comodule coalgebra. We recall and discuss the two notions in the first sections of the paper. The right action of $A$ on $B$ gives rise to the smash product $A\# B$. The left coaction of $B$ on $A$ gives a possible coproduct $\Delta_\#$ on $A\# B$. We will discuss in detail the necessary compatibility conditions between the action and the coaction for $\Delta_\#$ to be a proper coproduct on $A\# B$. The result is again a regular multiplier Hopf algebra. Majid's construction is obtained when we have Hopf algebras.
\snl
We also look at the dual case, constructed from a pair $(C,D)$ of regular multiplier Hopf algebras where now $C$ is a left 
$D$-module algebra while $D$ is a right $C$-comodule coalgebra. We will show that indeed, these two constructions are dual to each other in the sense that a natural pairing of $A$ with $C$ and of $B$ with $D$ will yield a duality between $A\# B$ and the smash product $C\# D$.
\snl
We show that the bicrossproduct of algebraic quantum groups is again an algebraic quantum group (i.e.\ a regular multiplier Hopf algebra with integrals). The $^*$-algebra case will also be considered. Some special cases will be treated and they will be related with other constructions available in the literature.
\snl
Finally, the basic example, coming from a (not necessarily finite) group $G$ with two subgroups $H$ and $K$ such that $G=KH$ and $H\cap K=\{e\}$ (where $e$ is the identity of $G$) will be used throughout the paper for motivation and illustration of the different notions and results. The cases where either $H$ or $K$ is a normal subgroup will get special attention.
\nl
\nl
March 2009 ({\it Version} 2.1)
\vskip 1 cm
\hrule
\bigskip\parindent 0.3 cm
\item{$^{\text{1}}$} Department of Mathematics, Hasselt University, Agoralaan, B-3590 Diepenbeek (Belgium). E-mail: Lydia.Delvaux\@uhasselt.be
\item{$^{\text{2}}$} Department of Mathematics, K.U.\ Leuven, Celestijnenlaan 200B (bus 2400),
B-3001 Heverlee (Belgium). E-mail: Alfons.VanDaele\@wis.kuleuven.be
\item{$^{\text{3}}$} Department of Mathematics, Southeast University, Nanjing 210096, China. E-mail: shuanhwang\@seu.edu.cn
\parindent 0 cm
 
\newpage

\bf 0. Introduction \rm
\nl
Consider a group $G$ and assume that it has two subgroups $H$ and $K$ such that $G=KH$ and $H\cap K=\{e\}$ where $e$ denotes the identity element in $G$. For any two elements $h\in H$ and $k\in K$, there is a unique way to write $hk$ as a product $k'h'$ with $h'\in H$ and $k'\in K$. We will write $h\tr k$ for $k'$ and $h\tl k$ for $h'$. As a consequence of the associativity of the product in $G$, we get, among other formulas, that
$$\align (hh')\tr k&= h\tr (h'\tr k) \\    
         h \tl (kk')&= (h\tl k)\tl k'
\endalign$$
for all $h,h'\in H$ and $k,k'\in K$. We get a left action of the group $H$ on the set $K$ and a right action of the group $K$ on the set $H$. It justifies the notations we have used. These actions are related in a very specific way. It is a (right-left) \it matched pair \rm $(K,H)$ of groups as defined by Majid (see e.g. Definition 6.2.10 in [M3]). We will use such a matched pair throughout the paper, mainly for motivational reasons. Therefore, we will recall more details about such a pair at various places in our paper.
\snl
Given a pair of groups $(K,H)$, we can associate two Hopf algebras. The first one is the group algebra ${\Bbb C}H$. We will use $A$ to denote this group algebra. The second one is the group algebra ${\Bbb C}K$. We will use $D$ to denote this group algebra. On the other hand, we can also consider the algebras $F(H)$ and $F(K)$ of complex functions with finite support on $H$ and $K$ respectively (with pointwise operations). We will use $B$ for $F(K)$ and $C$ for $F(H)$. If $H$ and $K$ are finite, also $B$ and $C$ are Hopf algebras and of course, $C$ is the dual of $A$ while $D$ is dual to $B$. In general, $B$ and $C$ are regular multiplier Hopf algebras with integrals, i.e.\ algebraic quantum groups (as defined and studied in [VD1] and [VD2]). Still, we have that $A$ and $C$ on the one hand and $B$ and $D$ on the other hand, are each others dual as algebraic quantum groups.
\snl
If moreover the pair $(K,H)$ is a matched pair of groups as above, the left action of $H$ on $K$ will induce a right action of $A$ on $B$ (in the sense of [Dr-VD-Z]). It is denoted and defined by
$$(f\tl h)(k)=f(h\tr k)$$
whenever $h\in H$, $k\in K$ and $f\in B$. This action makes $B$ into a right $A$-module algebra (see Section 1). Similarly, the right action of $K$ on $H$ will induce a left action of $D$ on $C$, denoted and defined by 
$$(k\tr f)(h)=f(h\tl k)$$
whenever $h\in H$, $k\in K$ and $f\in C$. It makes $C$ into a left $D$-module algebra.
\snl
As in [Dr-VD-Z], we can define the smash products $A\#B$ and $C\#D$. We will recall these notions in detail in Section 1 of this paper. In the case considered here in the introduction, the algebra $A\#B$ is the space $A\ot B$ with the product defined by
$$(h\ot f)(h'\ot f')=hh'\ot (f\tl h')f'$$
when $h,h'\in H$ and $f,f'\in B$. Observe that we consider elements of the group as sitting in the group algebra. Similarly, the algebra $C\#D$ here is the space $C\otimes D$  with product defined by
$$(f\ot k)(f'\ot k')=f(k\tr f')\ot kk'$$
whenever $k,k'\in K$ and $f,f'\in C$. See Example 1.8 in Section 1 for more details.
\nl
In the case of a pair of finite groups, the spaces are finite-dimensional. Then the multiplications on $A\ot B$  and $C\ot D$ induce comultiplications on $C\ot D$ and $A\ot B$ respectively by duality. It turns out that these coproducts make $A\#B$ and $C\#D$ into Hopf algebras, dual to each other. The compatibility of the actions of $H$ on $K$ and $K$ on $H$, as mentioned earlier in this introduction, is crucial for this (rather remarkable) result. This property was shown by Majid and a proof can be found in [M3] (see Example 6.2.11 and 6.2.12 in [M3]).
\snl
It is rather straightforward to show that this result of Majid remains valid when the groups are no longer assumed to be finite. Of course, in this more general case, we have to work with multiplier Hopf algebras as the spaces are no longer finite-dimensional so that the products on $A\ot B$ and $C\ot D$ will not induce ordinary coproducts on the dual spaces. However, apart from this, let us say technical problem, the extension of this result from the case of finite groups to the general case is, as said, mostly straightforward. This case has been considered in [VD-W]. Here, we will conclude it from general results on bicrossproducts for multiplier Hopf algebras, as obtained in the paper.
\nl
Indeed, the notion of a matched pair of groups is the underlying idea of Majid's bicrossproduct construction. For his construction, the starting point is a pair $(A,B)$ of Hopf algebras with a right action of $A$ on $B$ and a left coaction of $B$ on $A$. Roughly speaking, such a left coaction is the dual concept of a left action of $D$ (the dual of $B$) on $C$ (the dual of $A$) (see Section 2 in this paper). Again, there are compatibility conditions. These allow to show that the smash product $A\#B$ carries a coproduct (induced by the coaction), making it into a Hopf algebra.
\snl
The {\it main purpose of this paper} is to extend Majid's bicrossproduct construction to the case of multiplier Hopf algebras. We refer to the summary of the content of this paper, further in this introduction, for a more detailed account of what we do.
\snl
Still, the notion of a matched pair of groups (now possibly infinite), will be used as a motivation throughout our paper. And one might expect that the passage from Hopf algebras to regular multiplier Hopf algebras is again straightforward (as it is for the passage above from finite to infinite groups). To a certain extent, this is correct as essentially the formulas are the same. However, there are some difficulties as a consequence of the fact that, in the case of multiplier Hopf algebras, the coproducts (as well as the coactions) no longer take values in the tensor product of the algebras but rather in the multiplier algebras (see further in the introduction where we recall this notion very briefly).
\snl
There are essentially {\it two ways} to deal with this problem. One way is to keep using the same formulas, work with the Sweedler notations, both for the coproducts and the coactions, and verify carefully if everything is 'well-covered'. The {\it technique of covering} the legs of coproducts with values in the multiplier algebra (together with the use of Sweedler's notation) was first introduced in [Dr-VD]. However, in this context, we encounter a greater complexity with the use of the Sweedler notation. In [VD4] we have developed a more fundamental way to deal with these problems and we refer the reader to this article for a better theoretical basis for all this. In fact, in that paper, we have various examples, related to this work on bicrossproducts.
\snl
Certainly in this setting, the covering technique requires some care, but it has the advantage of being more transparent. Another possible way to study these bicrossproducts is by using {\it linear maps} between tensor products of spaces. It is like treating the coproduct $\Delta$ for a multiplier Hopf algebra $(A,\Delta)$ by using the associated linear maps $T_1, T_2$ from $A\ot A$ to $A\ot A$ given by
$$\align T_1(a\ot a') &= \Delta(a)(1\ot a') \\
         T_2(a\ot a') &= (a\ot 1)\Delta(a').  
\endalign$$
\snl
The two approaches look very different, but are in fact two forms of the same underlying idea. Roughly speaking one could say that the covering technique is a more transparent way of treating the correct formulas involving linear operators.
\snl
For all these reasons, we have chosen to write this paper with the following underlying point of view. We will spend time on motivation and doing so, we will focus on the use of the covering technique. On the other hand, we will also indicate how, in certain cases, this is all translated into rigorous formulas using linear maps. We do get important new results because we work in a different and more general setting, but the results are expected. By focusing more on the new techniques, we hope to make the paper also more interesting to read. Moreover, with this point of view, we move our paper slightly into the direction of a paper with a more expository nature. We feel this is also important because it can help the reader to get more familiar with this generalization of Hopf algebras which, after all, turns out to give nice and more general results. The bicrossproduct construction of Majid, treated in this paper, seems to be a suitable subject for this purpose.
\snl
Before we come to the content of the paper, we should mention in passing that Majid's bicrossproduct construction has also been obtained within the theory of locally compact quantum groups. The final result is obtained by Vaes and Vainerman [V-V2], see also [V1], [V2] and [V-V1]. Their work is greatly inspired by the paper by Baaj and Skandalis on multiplicative unitaries [B-S]. But before all this, we have the work of Majid on bicrossproducts for Hopf von Neumann algebras and Kac algebras [M1].
\snl
It is worthwhile mentioning that the theory of multiplier Hopf algebras and algebraic quantum groups has been developed before the present theory of locally compact quantum groups and that it has served as a source of inspiration for the latter. So, it should not come as a surprise that the present work is related with the earlier and more recent work on bicrossproducts in Hopf-von Neumann algebra theory. In fact, the formulas used in these analytical versions of the bicrossproduct construction are essentially the same as the ones that we use in this paper when we treat the matter with the 'linear operator technique'. We will say more about this at the appropriate place in the paper. Note however that our case can not be seen as a special case of this analytical theory because we work in a purely algebraic context. We will also come back to this statement later. See some remarks in Section 4.
\snl
We also should mention here that many interesting cases of bicrossproducts in the setting of locally compact quantum groups (as treated by Vaes and Vainerman in [V-V2]) do not fit into this algebraic approach. Nevertheless, there are interesting examples that cannot be treated using only  Hopf algebras but where multiplier Hopf algebras are needed (and sufficient). We will say something more in Section 5 of this paper. We also refer to our second paper on this subject [De-VD-W], where we include more examples. And observe also that the theory of multiplier Hopf algebras and algebraic quantum groups is, from a didactical point of view, a good intermediate step towards the more general, but also technically far more complicated theory of locally compact quantum groups. See e.g.\ reference [VD3] (Multiplier Hopf $^*$-algebras with positive integrals: A laboratory for locally compact quantum groups).
\nl
\it Content of the paper \rm
\nl
In {\it Section} 1, we recall briefly the notion of a unital action and the smash product. We first consider a right action of a regular multiplier Hopf algebra $A$ on an algebra $B$ and we assume that $B$ is a right $A$-module algebra. We recall the twist map (or braiding) and the notion of the smash product. The case of a left module algebra is also discussed. We explain a basic procedure to pass from one case to the other and to find the formulas in the case of a left module algebra from those of a right module algebra. We look at the $^*$-algebra cases as well. 
\snl
For the sake of completeness, because we will treat these also in the coming sections, we consider the trivial cases (when the actions are trivial). The example of two subgroups of a group (as we described earlier in this introduction) is the basic illustration of the notions. And also for this example, we look at the special case where one of the subgroups is normal.
\snl
We use this well-known and relatively simple situation to explain more precisely how to deal with problems that arise from the fact that we work with multiplier Hopf algebras where the coproduct has its range in the multiplier algebra of the tensor product and not in the tensor product itself. In particular, we recall how the Sweedler notation is used in this context and how to work with the {\it covering technique} properly. In [VD4], we have explained the theoretical background for these methods.
\snl
Actions of multiplier Hopf algebras and module algebras were introduced in [Dr-VD-Z]. Also the smash coproduct has been considered in several papers before (see also [De1] and [De3]). This first section does not contain really new ideas and certainly no new results. It is just included for convenience of the reader and as we mentioned, to illustrate the use of the Sweedler notation and the covering technique (as treated in the [VD4]).
\snl
In {\it Section} 2, we recall definitions and results about coactions and smash coproducts. The basic ingredient is that of a left coaction of a regular multiplier Hopf algebra $B$ on a vector space $A$. If $A$ is also an algebra, it is an injective linear map $\Gamma:A\to M(B\ot A)$, satisfying certain properties. In general, one must be a little more careful. One of the problems arise because $\Gamma$ has its range in the multiplier algebra of $B\ot A$ and not in $B\ot A$ itself. This causes the need for a new type of covering. Further, the associated cotwist map (or cobraiding) $T$ on $A\ot B$ is introduced (when also $A$ is a regular multiplier Hopf algebra). Formally, it gives rise to the coproduct $\Delta$ on $A\ot B$ given by $(\iota_A\ot T\ot \iota_B)(\Delta_A \ot \Delta_B)$. Again, there is a complication because the coproduct is supposed to map to the multiplier algebra of the tensor product and therefore, the definition of this coproduct requires, in some sense, already the product. Also this problem is discussed. 
\snl
The notion of a coaction was introduced in [VD-Z2] and in [De2] it is studied further (see also [De4]). In [De2], the cotwist map is considered and the notion of a comodule coalgebra is given. In this paper, we recall (and slightly generalize) these various notions and results. We also add some observations and provide some deeper insight in the problems that make this 'dual case' much more involved than the easier case of an action and a smash product as in Section 1. 
\snl
We consider the case of a right coaction as well. Again, we can use the basic technique as explained already in Section 1 to pass from one case to the other. Now, there is also a second possibility to do this, based on duality. This is also explained in this section. The two techniques are used further in Sections 3 and 4 to obtain the right formulas in one case from those in the other case. Also here, we treat the $^*$-algebra cases.
\snl
Finally in this section, again we consider the obvious trivial cases and we illustrate all of this with the group example as introduced before.
\nl
{\it Section} 3 is {\it the most important} section of the paper. In this section, we consider the candidate for the smash coproduct $\Delta_\#$ on the smash product $A\# B$. We explain step by step what the natural conditions are for $\Delta_\#$ to be an algebra map. The conditions involve different connections between the right action of $A$ on $B$ and the left coaction of $B$ on $A$. We formulate these conditions in a beautiful and symmetrical way, different from the classical formulations in Hopf algebra theory. And of course, we prove that under these conditions, we get indeed that $\Delta_\#$ is a coproduct on the smash product.
\snl
The other case of a left action of $D$ on $C$ and a right coaction of $C$ on $D$ is obtained as before, using the techniques we described in Section 1 and Section 2.
\snl
The $^*$-algebra case is treated, as well as the various special cases. The basic example is used to illustrate the conditions and the results.
\nl
In {\it Section} 4, we obtain the {\it main results} of this paper. We show that the smash product $A\# B$, as reviewed in the first section, endowed with the smash coproduct $\Delta_\#$, as defined in the second section, when the conditions discussed in the third section are fulfilled, is actually a regular multiplier Hopf algebra. Again also the dual case, the $^*$-algebra case and the case of algebraic quantum groups are considered, as well as the obvious special cases. In this section, the treatment of the basic example is completed. We also refer to a forthcoming paper where we complete the case of algebraic quantum groups in the sense that we also obtain formulas for the various objects like the modular elements, the modular automorphism groups, etc., associated with the bicrossproduct and its dual. See [De-VD-W].
\nl
In the last section, {\it Section 5}, we conclude the paper and we discuss some possible further research on this subject. In particular, we consider some remaining difficulties related with aspects of covering for coactions. And we discuss the problem of finding more examples that fit into this framework. 
\nl
We do not treat many examples in this paper. On the other hand, in our second paper on the subject [De-VD-W], where we consider integrals on bicrossproducts, we do consider various non-trivial and interesting examples.

\nl
\it Notations, conventions and basic references\rm
\nl
Throughout the paper, we work with (associative) algebras over the complex numbers $\Bbb C$ as often we are also interested in the $^*$-algebra case. We do not require the algebras to have an identity, but we do expect that the product is non-degenerate (as a bilinear map). For an algebra $A$, we use $M(A)$ to denote the multiplier algebra. It is the largest unital algebra, containing $A$ as a dense two-sided ideal. If $A$ is a $^*$-algebra, then so is $M(A)$.
\snl
If $A$ and $B$ are algebras and if $\alpha:A\to M(B)$ is an algebra homomorphism, then it is called non-degenerate if $\alpha(A)B=B$ and $B\alpha(A)=B$. In that case, $\alpha$ has a unique extension to a unital homomorphism from $M(A)$ to $M(B)$. This extension is still denoted by $\alpha$.
\snl
We use $1_A$ for the identity in $M(A)$ and simply $1$ if no confusion is possible. Sometimes however, we will even then use $1_A$ for clarity. The identity element in a group is denoted by $e$. We use $\iota_A$ for the identity map from $A$ to itself and again, we simply write $\iota$ when appropriate. Similarly, we use $\Delta_A$ and $\Delta$ for a coproduct on $A$.
\snl
A multiplier Hopf algebra is an algebra $A$ with a coproduct $\Delta$, satisfying certain assumptions. The coproduct is a non-degenerate homomorphism from $A$ to $M(A\ot A)$. It is assumed to be coassociative, i.e.\ we have $(\Delta\ot \iota)\Delta=(\iota\ot\Delta)\Delta$. If $A$ is a $^*$-algebra, we require that $\Delta$ is a $^*$-homomorphism. A multiplier Hopf algebra is called regular if the opposite coproduct $\Delta^{\text{op}}$ still makes $A$ into a multiplier Hopf algebra. For a multiplier Hopf $^*$-algebra, this is automatic. A regular multiplier Hopf algebra that carries integrals is called an algebraic quantum group. We refer to [VD1] for the theory of multiplier Hopf algebras and to [VD2] for the theory of algebraic quantum groups, see also [VD-Z1]. For the use of the Sweedler notation, as introduced in [Dr-VD], we refer to special paper about this subject, see [VD4]. For pairings of multiplier Hopf algebras, the main reference is also [Dr-VD]. More (basic) references will be given further in the paper. 

\nl\nl
\bf Acknowledgements \rm
\nl
We would like to thank the referee who reviewed an earlier version of this paper for many valuable remarks. They have inspired us to rewrite the paper more fundamentally. 
\snl
This work was partially supported by  
the FNS of CHINA (10871042). The third author also thanks the K.U.\ Leuven for the research fellowships he received in 2006 and 2009.
\nl

\bf 1. Actions and smash products \rm
\nl
In this section, we consider a regular multiplier Hopf algebra $A$ and any other algebra $B$. Recall the following notion (see [Dr-VD] and [Dr-VD-Z]).

\inspr{1.1} Definition \rm
Suppose that we have a unital right action of $A$ on $B$, denoted as $B\tl A$. Then $B$ is called a \it right $A$-module algebra \rm if also
$$bb'\tl a = \sum_{(a)}(b\tl a_{(1)})(b'\tl a_{(2)})$$
for all $a\in A$ and $b, b'\in B$. \hfill$\square$
\einspr

This notion has been studied before and is now well understood. However, because it is one of the basic ingredients of the construction in this paper and because of the scope of this paper, as outlined in the introduction, we will here discuss some aspects of this notion. In particular, we will focus on the technique of {\it covering}. Being a relatively simple concept, the notion of a module algebra is well adapted to explain something more about the covering technique. It should enable the reader to become more familiar with it.
\snl
That $A$ acts from the right on $B$ means in the first place that $b\tl (aa')=(b\tl a)\tl a'$ for all $a,a'\in A$. The action is assumed unital. This means that elements of the form $b\tl a$ with $a\in A$ and $b\in B$ span all of $B$. Because in any regular multiplier Hopf algebra, there exist local units, it follows that for any $b\in B$, there exists an element $e\in A$ such that $b\tl e=b$. Then, in the defining formula of the above definition, $a_{(1)}$ will be covered by $b$ and $a_{(2)}$ will be covered by $b'$. In Example 2.8.ii of [VD4], this case is used to illustrate the basic ideas behind the covering technique. 
\snl
As we also mentioned in the introduction, one way to avoid these coverings is by using linear maps on tensor product spaces. We will now again do this in greater detail for the notion of a module algebra as in Definition 1.1. For this purpose, a {\it twist map} is associated with the action. We recall the definition and the first main properties in the following proposition (see e.g.\ [Dr-VD-Z] where the case of a left module algebra is treated).

\inspr{1.2} Proposition \rm
Assume  that $B$ is a right $A$-module algebra. Define a linear map $R:B\ot A \to A\ot B$ by 
$$R(b\ot a)=\sum_{(a)}a_{(1)} \ot (b\tl a_{(2)}).$$
Then $R$ is bijective and 
$$R^{-1}(a\ot b)=\sum_{(a)}(b\tl S^{-1}(a_{(2)}))\ot a_{(1)}$$
for all $a\in A$ and $b\in B$. Here, $S$ is the antipode on $A$. \hfill$\square$
\einspr

Remark that in both equations, $a_{(2)}$ is covered by $b$ through the action. Again, the reader can look at Example 2.9.i in [VD4] for more details if desirable.
\snl
Some authors use a different terminology. In [B-M] e.g.\ this map is called a 'generalized braiding'. We prefer our terminology because the map is used to twist the usual product on $A\ot B$ (see Proposition 1.4 below). Also, it does not really satisfy the standard braid relation but rather the similar relation in Proposition 1.3 below.  
\snl
Indeed, using this linear map $R$, it is possible to translate the basic conditions of Definition 1.1 in terms of linear maps. These equations are given in the following proposition. See also e.g.\ [VD-VK] where these equalities appear in a natural way.

\inspr{1.3} Proposition \rm
We have
$$\align R(\iota_B\ot m_A)
      &=(m_A\ot \iota_B)(\iota_A\ot R)(R\ot \iota_A) \qquad \text{on} \qquad B\ot A\ot A \\
   R(m_B\ot\iota_A) 
      &=(\iota_A\ot m_B)(R\ot\iota_B)(\iota_B\ot R) \qquad \text{on} \qquad B\ot B\ot A
\endalign$$
where as mentioned before, $\iota_A, \iota_B$ are the identity maps on $A$ and on $B$ resp. and $m_A$ and $m_B$ are the multiplication maps. \hfill $\square$
\einspr

The first formula follows from the fact that $B$ is an $A$-module. The second one comes from the $A$-module algebra condition. This is easy to show (see e.g.\ [De3]). It can also be shown that these two equations in turn will imply that $B$ is a right $A$-module algebra. If e.g.\ we apply the first equation on $b\ot a\ot a'$ (with $a,a'\in A$ and $b\in B$) and then apply $\varepsilon_A\ot \iota_B$ we will obtain that $b\tl (aa')=(b\tl a)\tl a'$. Similarly we can get the other requirement from the second equation. 

\nl
Next we consider the \it smash product\rm. There are different ways to treat this concept. We start with recalling the usual approach (see [Dr-VD-Z] for the case of a left module algebra and [De3]). Of course, also elsewhere in the literature, related work has been done.

\inspr{1.4} Proposition \rm
Define a linear map $m:(A\ot B)\ot (A\ot B) \to A\ot B$ by
$$m=(m_A\ot m_B)(\iota_A\ot R \ot \iota_B).$$
This map makes $A\ot B$ into an (associative) algebra (with a non-degenerate product).

\hfill$\square$
\einspr

The associativity can be proven by using the two formulas for $R$ obtained in Proposition 1.3. 
\snl
First, we will follow the common convention and write $A\# B$ for the space $A\ot B$, endowed with this product. Similarly, we will use $a\# b$ for the element $a\ot b$ when considered as an element of this algebra.
\snl
Using Sweedler's notation, the product can be written as
$$(a\#b)(a'\#b')=\sum_{(a')}a a'_{(1)}\# (b\tl a'_{(2)})b'$$
for all $a,a'\in A$ and $b,b'\in B$. In this formula, $a'_{(1)}$ is covered by $a$ (through multiplication) and $a'_{(2)}$ is covered by $b$ (through the action).
\snl
We are also interested in the other way to treat this smash product because it makes formulas sometimes more transparent. It is based on the following result (see Proposition 5.10 in [Dr-VD-Z]).

\inspr{1.5} Proposition \rm 
There exist injective non-degenerate homomorphisms \newline
$\pi_A:A\to M(A\# B)$ and $\pi_B:B\to M(A\# B)$ defined by
$$\align \pi_A(a')(a\# b)&=(aa')\# b \\
         (a\# b)\pi_B(b')&=a\# (bb')
\endalign$$
for all $a,a'\in A$ and $b,b' \in B$. We also have
$$\align \pi_A(a)\pi_B(b)&=a\# b \\
         \pi_B(b)\pi_A(a)&=\sum_{(a)}a_{(1)}\# (b\tl a_{(2)})
\endalign$$
for all $a\in A$ and $b\in B$. \hfill $\square$
\einspr
Observe the presence of the map $R$ in the right hand side of the last formula.
\snl
If we would identify $A$ and $B$ with their images in $M(A\# B)$, we can replace $a\# b$ by $ab$ and we find that $ba=\sum_{(a)}a_{(1)} (b\tl a_{(2)})$ for all $a\in A$ and $b\in b$. This takes us to the following result that provides an alternative approach to the smash product.

\inspr{1.6} Proposition \rm 
The smash product algebra $A\# B$ is isomorphic with the algebra, generated by $A$ and $B$, subject to the commutation rules $ba=\sum_{(a)}a_{(1)}(b\tl a_{(2)})$ for all $a\in A$ and $b\in B$. \hfill$\square$
\einspr

It makes sense to denote this algebra with $AB$ and it requires a small argument (bijectivity of the map $R$) to show that also $AB=BA$. More precisely, the two maps $a\ot b\mapsto ab$ and $b\ot a\mapsto ba$ are linear bijections from the spaces $A\ot B$ and $B\ot A$ respectively to the space $AB$. In the remaining part of the paper, we will use $AB$ more often than $A\# B$ to denote the smash product.
\snl
The algebra $AB$ acts faithfully on itself by right multiplication. A simpler action is obtained on $B$ by letting $A$ act in the original way and $B$ again by right multiplication. Also this action will be used further in the paper although in general, it need not be faithful anymore (e.g.\ when the action of $A$ on $B$ is trivial).
\nl
If the action is trivial, that is if $b\tl a=\varepsilon(a)b$ for all $a\in A$ and $b\in B$, then the twist map is nothing else but the flip and the smash product is simply the tensor product algebra. In the other approach, it means that $A$ and $B$ commute in the smash product $AB$.
\nl
If $A$ is a multiplier Hopf $^*$-algebra and if $B$ is a $^*$-algebra, we require the extra condition that
$$(b\tl a)^*=b^*\tl {S(a)}^*$$
for all $a\in A$  and $b\in B$. This is a natural condition. One easily verifies e.g.\ that it is compatible, both with the action property (as ${S(aa')}^*=S(a)^*{S(a')}^*$ for all $a,a'\in A$), and the module algebra property (as $(bb')^*={b'}^*b^*$ for all $b,b'\in B$ and $\Delta(S(a)^*)=((S\ot S)\Delta^{\text{op}}(a))^*$ for all $a\in A$). One can show that with this extra assumption, the smash product $A\# B$ can be made into a $^*$-algebra simply by letting $(ab)^*=b^*a^*$ for all $a\in A$ and $b\in B$.
\nl
Before we consider the basic example, we also look (very briefly) at the case of a left $D$-module algebra $C$. We have made it clear in the introduction why we also need this case. Definitions and results can be taken directly from [Dr-VD-Z].

\inspr{1.7} Definition  \rm 
Let $D$ be a regular multiplier Hopf algebra and $C$ an algebra. Assume that $D$ acts from the left on $C$ (denoted $D\tr C$). Then $C$ is called a {\it left $D$-module algebra} if also 
$$d\tr (cc')=\sum_{(d)}(d_{(1)}\tr c)(d_{(2)}\tr c')$$
for all $c,c'\in C$ and $d\in D$. \hfill$\square$
\einspr

The relevant twist map $R$ is now a linear map from $D\ot C$ to $C \ot D$ and it is given by
$$R(d\ot c)=\sum_{(d)}(d_{(1)}\tr c)\ot d_{(2)}.$$
The smash product $C\# D$ is the algebra generated by $C$ and $D$, subject to the commutation rules 
$$dc=\sum_{(d)}(d_{(1)}\tr c)d_{(2)}$$
for all $c\in C$ and $d\in D$. It is denoted by $CD$. Again the maps $c\ot d\mapsto cd$ and $d\ot c \mapsto dc$ are bijective
from the spaces $C\ot D$ and $D\ot C$ respectively to the space $CD$.
\snl
Remark that there is some possible confusion because we use the same notations as in the case of a right $A$-module algebra $B$. A similar problem will occur on later occasions. However, we will be consequent in the use of $A$ and $B$ in the first case and $C$  and $D$ in the second case. This should help to distinguish the cases.
\snl
As for a right module algebra, also here, we have the $^*$-algebra case. Then we assume that $D$ is a multiplier Hopf $^*$-algebra, that $C$ is a $^*$-algebra and we make the extra assumption that $(d\tr c)^*=S(d)^*\tr c^*$ for all $c\in C$ and $d\in D$. Now the smash product $C\# D$ is again a $^*$-algebra with the involution defined by $(cd)^*=d^*c^*$
\nl
There is a procedure to convert formulas from one case, say the right module algebras, to the other, now the left module algebras. Indeed, if $B$ is a right $A$-module algebra (as in Definition 1.1) then we obtain a left $D$-module algebra $C$ if we let $C$ be the algebra $B$ endowed with the opposite product and if we take for the multiplier Hopf algebra $D$ the algebra $A$ but with the opposite product and coproduct. If we also flip the notation for the action, as well as the tensor products, we see that the formula in Definition 1.1 becomes the one in Definition 1.7, that the formula for $R$ in Proposition 1.2 converts to the formula for $R$ as above in the other case, etc. Remark however that the formulas in Proposition 1.3 look the same, but that they are in fact interchanged.
\snl
This is an important technical tool for this paper because we need to work with both cases. 
\nl
Now we consider the main motivating example, as announced already in the introduction.

\inspr{1.8} Example \rm 
Let $H$ and $K$ be groups.
\snl
i) Assume first that $H$ acts on $K$ from the left (as a group on a set). We use $h\tr k$ to denote the action of an element $h\in H$ on a point $k\in K$. We assume that the action is unital in the sense that $e\tr k=k$ for all $k\in K$ where $e$ denotes the unit in the group $H$. Let $A$ be the group algebra $\Bbb C H$ of $H$. It is a Hopf algebra (and in particular, a regular multiplier Hopf algebra). Let $B$ denote the algebra $F(K)$ of complex functions on $K$ with finite support, endowed with pointwise operations. The left action of $H$ on $K$ induces a right action of $A$ on $B$ by
$$(f\tl h)(k)=f(h\tr k)$$
whenever $h\in H$, $k\in K$ and $f\in B$. Remark that we consider elements of the group $H$ as sitting in the group algebra $\Bbb C H$ as usual. Then $B$ is a right $A$-module algebra in the sense of Definition 1.1.
\snl
The smash product  $AB$ is spanned by elements of the form $hf$ with $h\in H$ and $f\in B$ and the product is given by
$$(hf)(h'f')=(hh')((f\tl h')f')$$ 
whenever $h,h'\in H$ and $f,f'\in B$.
\snl
ii) Now assume that $K$ acts on $H$ from the right and use $h\tl k$ to denote the action of $k\in K$ on $h\in H$. Again assume that the action is unital. Let $C$ be the function algebra $F(H)$ and let $D$ be the group algebra $\Bbb C K$. The right action of $K$ on $H$ gives a left action of $D$ on $C$ defined by
$$(k\tr f)h=f(h\tl k)$$
for all $h\in H$, $k\in K$ and $f\in C$, making $C$ into a left $D$-module algebra. In this case, the smash product $CD$ is spanned by elements of the form $fk$ with $f\in C$ and $k\in K$ and the product is given by 
$$(fk)(f'k')=(f(k\tr f'))(kk')$$
for all $k,k'\in K$ and $f,f'\in C$. \hfill$\square$
\einspr
Observe that in the above examples, the algebras $A$ and $D$ are Hopf algebras and we do not need to think about coverings.  
\snl
Also, when we endow these algebras with the obvious $^*$-algebra structures, the actions satisfy the requirements. We have e.g.\ that $f^*=f^-$, the complex conjugate of $f$ when $f\in F(K)$ and $h^*=h^{-1}$ so that $S(h)^*=h$ for $h\in H$. This will give the equality $(f\tl h)^*=f^*\tl S(h)^*$. 
\snl
Finally, if the group actions are trivial on the sets, they also are trivial on the algebras (remark that the counit maps each group element to $1$). Then the smash products are simply the tensor products. 
\nl\nl

\bf 2. Coactions and smash coproducts \rm
\nl
In this section, we will consider the objects, dual to actions and smash products as reviewed in the previous section. We will state and discuss notions and results as they  have appeared elsewhere in the literature (see [VD-Z2] and [De2]). However, we start the treatment in a more general (and perhaps also more natural) setting.  Further, we focus on motivation and on those aspects that are important for the main construction of the bicrossproduct in the next two sections. In particular, as we also have done in the previous section, we will be concerned with the proper coverings of our formulas and expressions therein.
\snl
As we will see, the situation is somewhat more involved than for actions and smash products.
\nl
\it Coactions \rm
\nl
It is natural to start in this case with the notion of a coaction. It has been introduced in [VD-Z2] and studied further in [De2]. In this paper, we will treat a slightly more general case. In Hopf algebra theory, it is possible to define a coaction of a coalgebra on a vector space. In the setting of multiplier Hopf algebras however, more structure is needed. In [VD-Z2], the setting is that of an algebra $A$ and a regular multiplier Hopf algebra $B$. Then, a left coaction of $B$ on $A$ is an injective linear map $\Gamma:A\to M(B\ot A)$ satisfying 
\snl 
\qquad i) \ $\Gamma(A)(B\ot 1)\subseteq B\ot A$ \quad and \quad $(B\ot 1)\Gamma(A) \subseteq B\ot A,$ 
\snl
\qquad ii) \ $(\iota_B\ot \Gamma)\Gamma = (\Delta_B\ot \iota_A)\Gamma.$
\snl
The algebra structures of $A$ and $B$ are needed to be able to consider the multiplier algebra $M(B\ot A)$. It would be too restrictive to assume that the coaction $\Gamma$ has range in the tensor product itself (see e.g.\ the example at the end of this section). Condition i) is used to give a meaning to the left hand side of the equation in the second condition. And to give a meaning to the right hand side of the second equation, the non-degeneracy of $\Delta_B\ot \iota_A$ is used in order to extend it to $M(B\ot A)$.
\snl
For the new notion here, we begin with the following notation. We refer to Section 1 in [VD4] for more details. In particular, see Definition 1.5 and Proposition 1.6 in [VD4] for the notion of an extended module of a bimodule.

\inspr{2.1} Notation \rm 
Let $A$ be a vector space and let $B$ be an algebra. Consider the $B$-bimodule $B\ot A$ with $B$ acting as multiplication left and right on the first factor in the tensor product. Denote by $M_0(B\ot A)$ the {\it completion of this module}. Similarly, consider $B\ot B\ot A$ as a $(B\ot B)$-bimodule with multiplication left and right on the first two factors. Denote by $M_0(B\ot B\ot A)$ the completion of this module. \hfill$\square$
\einspr

We will write the actions of $B$ as $(b\ot 1)y$ and $y(b\ot 1)$ when $b\in B$ and $y$ is an element of $B\ot A$ or $M_0(B\ot A)$. By definition, these elements belong to $B\ot A$, also for $y\in M_0(B\ot A)$. Similarly, we write the actions of $B\ot B$ as $(b\ot b'\ot 1)z$ and $z(b\ot b'\ot 1)$ when $b, b'\in B$ and $z$ is an element of either $B\ot B\ot A$ or $M_0(B\ot B\ot A)$. Again by definition, we get elements in $B\ot B\ot A$, also for $z\in M_0(B\ot B\ot A)$. This convention is in agreement with the notions that are commonly used when both $A$ and $B$ are algebras so that $M_0(B\ot A)$ and $M_0(B\ot B\ot A)$ are natural subspaces of $M(B\ot A)$ and $M(B\ot B\ot A)$ respectively.

\snl
Then we work with the following definition.

\inspr{2.2} Definition \rm
Let $A$ be a vector space and $B$ a regular multiplier Hopf algebra. A {\it left coaction} of $B$ on $A$ is an injective linear map $\Gamma:A\to M_0(B\ot A)$ such that $(\iota_B\ot \Gamma)\Gamma=(\Delta_B\ot \iota_A)\Gamma$ on $A$. We call $A$ a {\it left $B$-comodule}. \hfill$\square$
\einspr

It is easy to show that the map $\iota_B\ot\Gamma$ has a natural extension to $M_0(B\ot A)$. We simply multiply with an element of $B$ in the first factor. Also $\Delta_B\ot\iota_A$ can be extended to $M_0(B\ot A)$ in a natural way. Now multiply with $b\ot b'$ left or right. In the first case, use that $b\ot b'\in (B\ot 1)\Delta(B)$ and in the other case that $b\ot b'\in \Delta(B)(B\ot 1)$. These extensions will map into $M_0(B\ot B\ot A)$ and the condition $(\iota_B\ot \Gamma)\Gamma(a)=(\Delta_B\ot \iota_A)\Gamma(a)$ is a well-defined equation in this space. For more details, again we refer to Section 1 in [VD4].
 
\snl
It does not seem natural to weaken the assumptions about $B$ (as we can see above: we use that it is a regular multiplier Hopf algebra). On the other hand, remark that in the case where $A$ is an algebra, the notion coincides with the one originally given in [VD-Z2] and recalled in the beginning of this section. It is clear that the assumption i) in the original definition above has become part of the definition by the use of the space $M_0(B\ot A)$.
\snl
It is an easy consequence of the assumptions that $(\varepsilon_B\ot \iota_A)\Gamma(a)=a$ for all $a\in A$. The formula is given a meaning in the obvious way and the proof is the same as in the more restricted situation (see [VD-Z2] and [De2]). See also Example 2.10.i in [VD4].
\snl
We will also use the adapted version of the {\it Sweedler notation} for such a coaction. We will write
$\Gamma(a)=\sum_{(a)} a_{(-1)}\ot a_{(0)}$ so that coassociativity is implicit in the formula $(\iota_B\ot \Gamma)\Gamma(a)=
\sum_{(a)}a_{(-2)}\ot a_{(-1)} \ot a_{(0)}$ when it is understood that $\Delta_B(a_{(-1)})$ is replaced by $a_{(-2)}\ot a_{(-1)}$. We have more or less the same rules for covering (due to the two inclusions in i) above or, in the more general case, because of the definition). We need to cover the factor $a_{(-1)}$ and possibly also $a_{(-2)}$ (and so on) by elements in $B$, left or right. In this case however, one can not cover the factor $a_{(0)}$.
\snl 
One might expect that now, we can immediately move to the definition of a comodule coalgebra. Remember that in the previous section, there was no problem to define a module algebra already in the first definition. Here however, things are again more complicated.
\nl
\it The cotwist map $T$ \rm
\nl
First, we will make a remark that explains why we can consider coactions as dual to actions. This point of view will be used several times for motivation and later, in Sections 3 and 4, it will really be used to prove results about duality. We will see how this gives the natural candidate for the cotwist map.

\inspr{2.3} Remark \rm
Assume that we have regular multiplier Hopf algebras $A, B, C$ and $D$ and that $A$ is paired with $C$ and that $B$ is paired with $D$ (in the sense of [Dr-VD]). Use $\langle\,\cdot\, , \,\cdot\,\rangle$ to denote these pairings. Assume that we have a left action of $D$ on $C$. 
\snl
If all spaces are finite-dimensional, we can define a map $\Gamma:A\to B\ot A$ by
$$\langle \Gamma(a),d\ot c\rangle = \langle a, d \tr c \rangle$$
for all $a\in A$, $c\in C$ and $d\in D$. It is easily verified that 
$(dd')\tr c= d\tr (d'\tr c)$
will be the same as the property
$(\iota_B\ot \Gamma)\Gamma=(\Delta_B\ot\iota_A)\Gamma$.
\snl
Let us also look at the map $R$, now defined from $D\ot C$ to $C\ot D$ by 
$$R(d\ot c)=\sum_{(d)} (d_{(1)}\tr c) \ot d_{(2)}.$$
If again we assume our spaces to be finite-dimensional, we can look 
at the adjoint map $T:A\ot B\to B\ot A$. We find that 
$$\align \langle a\ot b, R(d\ot c)\rangle &= \sum_{(d)}\langle a, d_{(1)}\tr c \rangle \langle b, d_{(2)}\rangle \\
    &=\sum_{(d)}\langle \Gamma(a), d_{(1)} \ot c\rangle \langle b, d_{(2)} \rangle \\
    &=\langle \Gamma(a)(b\ot 1), d\ot c\rangle
\endalign$$
for all $a\in A$, $b\in B$, $c\in C$ and $d\in D$ and so $T(a\ot b)=\Gamma(a)(b\ot 1)$ for all $a\in A$ and $b\in B$.
\snl
Later in the paper, we will also have this type of duality in the more general case of pairings between regular multiplier Hopf algebras (see Section 4). \hfill$\square$
\einspr


The above argument suggests to define the cotwist map $T$. It will indeed play the role of $R$ in this dual setting. We have the following result. It is just as in the more restrictive case, considered in [De2].

\inspr{2.4} Proposition \rm 
Consider a left coaction $\Gamma$ as in Definition 2.2. 
Define a linear map $T:A\ot B\to B\ot A$ by $T(a\ot b)=\Gamma(a)(b\ot 1)$. Then $T$ is bijective and $T^{-1}$ is given by
$$T^{-1}(b\ot a)=\sum_{(a)}a_{(0)}\ot S^{-1}(a_{(-1)})b$$
for $a\in A$ and $b\in B$ (and where we use the Sweedler notation for $\Gamma(a)$ as explained earlier).
Moreover, 
$$(\Delta_B\ot \iota_A)T=(\iota_B\ot T)(T\ot \iota_B)(\iota_A\ot \Delta_B)$$
on $A\ot B$. \hfill$\square$
\einspr

There is no problem with the interpretation of the right hand side of the first formula as $a_{(-1)}$ will be covered through $b$ (in fact by $S(b)$). To prove that indeed, we get the inverse of $T$, we need to show that e.g.\
$$\sum_{(a)}a_{(-1)}S^{-1}(a_{(-2)})b \ot a_{(0)}=\sum_{(a)}\varepsilon_B(a_{(-1)})b\ot a_{(0)}$$
for all $a\in A$ and $b\in B$. In order to do this properly, we need to multiply in the first factor with an element in $B$ from the left. Doing so, we have covered the first two factors in the expression $(\Delta_B\ot \iota_A)\Gamma(a)$, we get something in $B\ot B\ot A$ and we can replace the expression $a_{(-1)}S^{-1}(a_{(-2)})$ with $\varepsilon_B(a_{(-1)})$. See e.g.\ Example 2.10.ii in [VD4].
\snl
In the second equation, when applied to $a\ot b$, we obviously can cover the left hand side of the equation with $b'\ot 1\ot 1$. Also the right hand side will be covered if we multiply with this element from the right. We use that
$$((T\ot\iota_B)(a\ot \Delta_B(b)))(b'\ot 1\ot 1)=(T\ot \iota_B)((a\ot\Delta_B(b))(1\ot b'\ot 1)).$$
The equality is then easily obtained from the coassociativity rule in Definition 2.2.
\snl
This last formula, involving $T$ and $\Delta_B$, is completely expected when we think of $T$ as the adjoint of the map $R$ in the setting of duality. It is simply the adjoint of the equation
$$R(m_D\ot\iota_C) =(\iota_C\ot m_D)(R\ot\iota_D)(\iota_D\ot R) \quad \text{on} \quad D\ot D\ot C,$$
a formula which is the counterpart of the first formula in Proposition 1.3, but now for a left $D$-module $C$.
\snl
Here, as you might expect, this formula is essentially equivalent with the coassociativity for the coaction, just as the first formula in Proposition 1.3 is essentially equivalent with the module property. But observe again that the formula for $T$ is of greater complexity than the formula for $R$, simply because in the last case, we do not have to worry about coverings.
\snl
We call this map a cotwist (as it will be used to twist the coproduct). Sometimes, the name (generalized) cobraiding is used (see e.g.\ [B-M]).
\nl
\it The $B$-comodule coalgebra $A$ \rm
\nl
Now we look for the right notion of a comodule coalgebra. For this, we assume that also $A$ is a regular multiplier Hopf algebra. 
\snl
It would be most natural to require, just as in the case of Hopf algebras, that
$$(\iota_B\ot\Delta_A)\Gamma(a)=\sum_{(a)}\Gamma_{12}(a_{(1)})\Gamma_{13}(a_{(2)}),$$
where we use the {\it leg numbering notation}. So, in this case, for $p\in A$, we have $\Gamma_{12}(p)=\Gamma(p)\ot 1$ and $\Gamma_{13}(p)=(\iota_B\ot\sigma_A)(\Gamma(p)\ot 1)$ where $\sigma_A$ is the flip on $A\ot A$.
When using the Sweedler notation, this equation reads as
$$(\iota_B\ot\Delta_A)\Gamma(a)=\sum_{(a)} a_{(1)(-1)}a_{(2)(-1)}\ot a_{(1)(0)}\ot a_{(2)(0)}$$
for all $a\in A$. In this last expression, the Sweedler notation is used, first for $\Delta_A$ and then for $\Gamma$. However, there is a serious problem with the coverings here. We can multiply with an element $b$ in the first factor and this will cover the $\Gamma$'s. Still, it seems not obvious how to cover $a_{(1)}$ and $a_{(2)}$ in the right hand side of the equation at the next stage. Because of this problem, we will not look at details here. We will consider this matter later.  
\snl
Another, still more or less obvious choice (for defining the notion of a comodule coalgebra) is to look at the adjoint of the second formula in Proposition 1.3 (for a left action). Then we would require 
$$(\iota_B\ot \Delta_A)T=(T\ot\iota_A)(\iota_A\ot T)(\Delta_A\ot\iota_B) \quad \text{on}\quad A\ot B.$$
This means that we essentially have used the covering with the element in $B$ as before and so it is expected that the problem remains. This is indeed the case. There is no problem with the left hand side of the equation, but it is still not obvious how to treat the right hand side.
\snl
However, it is now more easy to see how this problem can be overcome. The solution will be found from the result in the following lemma. 

\inspr{2.5} Lemma \rm First consider $A\ot B\ot A$ as an $A$-bimodule where the action of $A$ is given by multiplication in the first factor. Denote the completed module by $M_0^1(A\ot B\ot A)$. Then, one can define the map $(\iota_A\ot T)(\Delta_A\ot \iota_B)$ from $A\ot B$ to $M_0^1(A\ot B\ot A)$ in a natural way. Next, consider $A\ot B\ot A$ as an $A$-bimodule with multiplication in the third factor and denote the completed module by $M_0^3(A\ot B\ot A)$. Then there is a natural map $(T^{-1}\ot\iota_A)(\iota_B\ot\Delta_A)$ from $B\ot A$ to $M_0^3(A\ot B\ot A)$.
\hfill $\square$
\einspr

In the first case, define e.g.\
$$(a'\ot 1\ot 1)((\iota_A\ot T)(\Delta_A(a)\ot b))=(\iota_A\ot T)((a'\ot 1)\Delta_A(a)\ot b)$$
for all $a,a'\in A$ and $b\in B$. Similarly on the other side. This proves the first statement. The second one is proven in the same way. 
\snl
Given this result, it is possible to impose the required condition under the following form. Indeed, it makes sense to assume that 
$$(\iota_A\ot T)(\Delta_A\ot \iota_B)=(T^{-1}\ot\iota_A)(\iota_B\ot\Delta_A)T$$
on $A\ot B$. The left hand side gives an element in $M_0^1(A\ot B\ot A)$ while the right hand side an element in $M_0^3(A\ot B\ot A)$. Because there is a natural intersection of these two spaces, from this equality, the following assumption would be an implicit consequence. Therefore, let us make this assumption:

\inspr{2.6} Assumption \rm We assume that the elements of the form 
$$\align & ((\iota_A\ot T)(\Delta_A(a)\ot b))(1\ot 1\ot a') \\
         & (1\ot 1 \ot a') (\iota_A\ot T)(\Delta_A(a)\ot b)
\endalign$$
are in $A\ot B\ot A$ for all $a,a'\in A$ and $b\in B$.  \hfill$\square$
\einspr

We see from the discussion before, that this is a natural condition. We will say more about this condition after we have stated the following definition (cf. Definition 1.6 in [De2]). 

\inspr{2.7} Definition \rm 
Let $A$ and $B$ be regular multiplier Hopf algebras and assume that $\Gamma:A\to M(B\ot A)$ is a left coaction of $B$ on $A$ (as in Definition 2.2). Then we call $A$ a {\it left $B$-comodule coalgebra} if also
$$(\iota_B\ot\Delta_A)T=(T\ot \iota_A)(\iota_A\ot T)(\Delta_A\ot\iota_B)$$
on $A\ot B$. \hfill$\square$
\einspr

Compare this equality with the one in Proposition 2.4. It is very similar, but there, the covering problem was more easily solved. It is however also possible to treat the covering of the formula in Proposition 2.4 in a similar way as we did for the formula in Definition 2.7. It will not be necessary to impose an extra condition as in Assumption 2.6. as the analogue of this assumption will be automatic.
\snl
Before we proceed, let us make a {\it few more remarks} about the Assumption 2.6. 
\snl
First, it is easy to reformulate the assumption in terms of the map $\Gamma$ itself. We simply have that the expressions
$$((\iota_A\ot \Gamma)\Delta_A(a))(1\ot b\ot 1),$$
when multiplied with elements of $A$ in the third factor, left or right, give a result in the tensor product $A\ot B\ot A$.
\snl
Next, fix $a'\in A$ and $b\in B$ and look at the map $F:A\to B\ot A$ given by $F(a)=\Gamma(a)(b\ot a')$. If in this map, the variable $a$ is covered from the right, then the first part of the assumption is fulfilled. Similarly, take the map $G:A\to B\ot A$ given by $G(a)=(1\ot a')\Gamma(a)(b\ot 1)$ and assume that the variable $a$ is covered from the left (or from the right). Then, the second part of the assumption is fulfilled. 
\snl
There are reasons to believe that also the converse is true. Assume e.g.\ that $A$ is an algebraic quantum group of discrete type so that there is a left cointegral $h$. We know that the right leg of $\Delta(h)$ is all of $A$. Then, from the assumption, with $a=h$, it will follow that the maps $F$ and $G$ above have variables that are covered. As the result is trivially true for algebraic quantum groups of compact type, one might expect that it is true for all algebraic quantum groups, possibly for all regular multiplier Hopf algebras. 
\snl
Remark also that in the special case when $\Gamma$ is a homomorphism, then the existence of local units for $A$ together with the fact that $B\ot A=\Gamma(A)(B\ot 1)$, will imply that these variables are covered and hence that the Assumption 2.6 is automatic.
\snl
In fact, it would be convenient to assume the stronger assumption because then it would be possible to cover the right hand side in the original equations
$$(\iota_B\ot\Delta_A)\Gamma(a)=\sum_{(a)}\Gamma_{12}(a_{(1)})\Gamma_{13}(a_{(2)})
=\sum_{(a)} a_{(1)(-1)}a_{(2)(-1)}\ot a_{(1)(0)}\ot a_{(2)(0)}$$
by multiplying with $b\ot 1\ot a'$ where $a,a'\in A$ and $b\in B$. But for the moment, we see no good reason why we should assume this stronger property. We will just require Assumption 2.6. Then, if we multiply the above expressions with elements of the form $b\ot 1\ot a'$, from the right, where $b\in B$ and $a'\in A$, we do get everything well-defined in $B\ot A\ot A$. And although this does not really fit in the covering approach, this will be sufficient for our purposes. We will come back to this point in the next section where this will be used (see e.g.\ the proof of Proposition 3.11). We refer to Section 5 for some further discussion on this problem.
\nl
Before we proceed to the next item, just remark that it easily follows from the assumptions that $(\iota_B\ot\varepsilon_A)\Gamma(a)=\varepsilon_A(a) 1_B$ for all $a\in A$ (see e.g.\ Proposition 1.8 in [De2]). 
\nl
\it The twisted coproduct on $A\ot B$ \rm
\nl
The next logical step is the introduction of the {\it twisted coproduct}.
\snl
As explained before (in the introduction), in this context, we can not define a coproduct without referring to the algebra structure. This means that we expect that it is not possible to associate a twisted coproduct on $A\ot B$ when we (only) know that $A$ is left $B$-comodule coalgebra. The problem of course does not exist in the framework of Hopf algebras.
\snl
If we would consider the ordinary tensor product algebra structure on $A\ot B$, we can rely on the work done in [De2]. However, this is not the algebra structure on $A\ot B$ that we are interested in. Nevertheless, let us try and see how far we can get. We will therefore recall some of the results from [De2].
\snl
For, there is of course a {\it natural candidate} for the coproduct. Formally, we should define $\Delta$ on $A\ot B$ by the formula
$$\Delta(a\ot b)=(\iota_A\ot T \ot \iota_B)(\Delta_A(a) \ot \Delta_B(b)).$$
This is indeed dual to the formula of the product in the case of a module algebra (as given in Definition 1.4). Still formally, using Sweedler's notation, we can write
$$\Delta(a\ot b)=\sum_{(a)(b)}a_{(1)} \ot \Gamma(a_{(2)})(b_{(1)}\ot 1) \ot b_{(2)}$$
for $a\in A$  and $b\in B$.
\snl
The right hand side is easily covered when multiplying with elements of $A$ in the first factor, left or right, and with elements of $B$ in the last factor, also left or right. Unfortunately, this is not the type of covering that is needed to properly define a coproduct on the smash product. However, we also have the following possibility to cover this expression.

\inspr{2.8} Lemma \rm Take $A\ot B$ with the tensor product algebra structure. Consider $A\ot B\ot A\ot B$ as a $(A\ot B)$-bimodule where the actions are given by multiplication in the first two factors from the left and in the last two factors from the right. Then, there is a natural map $(\iota_A\ot T \ot \iota_B)(\Delta_A \ot \Delta_B)$ from $A \ot B$ to the completed module $M_0(A\ot B\ot A\ot B)$.

\snl \bf Proof: \rm 
First multiply from the left with such elements. For all $a,a'\in A$ and $b\in B$, we have
$$\align (a'\ot b'\ot 1 \ot 1)\sum_{(a)}(a_{(1)} \ot \Gamma(a_{(2)}) \ot 1) 
           &\subseteq (1\ot b'\ot 1 \ot 1)(A\ot \Gamma(A) \ot 1) \\
           &\subseteq A\ot B\ot A\ot 1.
\endalign$$
The second factor contains elements of $B$ and this will cover the factor $b_{(1)}$ in the expression
$$\Delta(a\ot b)=\sum_{(a)(b)}a_{(1)} \ot \Gamma(a_{(2)})(b_{(1)}\ot 1) \ot b_{(2)}.$$ 
\snl
Next, we multiply with an element of the form $1\ot 1\ot a'\ot b'$ from the right. Clearly, the element $b'$ will take care of the covering of the factor $b_{(2)}$. We get a linear combination of elements of the form
$$\sum_{(a)}a_{(1)} \ot \Gamma(a_{(2)})(p\ot a') \ot q$$
with $p,q\in B$. Then this belongs to $A\ot B\ot A\ot B$ by Assumption 2.6.   \hfill $\square$
\einspr

Similarly, we can multiply with elements of $A\ot B$ from the right in the first two factors and we will get elements in $A\ot B\ot A\ot B$. However, if we multiply with such elements from the left in the last two factors, we would need an assumption, similar as  Assumption 2.6, but for the opposite map $T^{\text{op}}$, defined by $T^{\text{op}}(a\ot b)=(b\ot 1)\Gamma(a)$. We will not need this for the moment.
\snl
It is an easy consequence of this result that the coproduct can be properly defined on $A\ot B$ with the above formula in the sense that we have a well-defined linear map $\Delta: A\ot B \to M((A\ot B)\ot (A\ot B))$ when $A\ot B$ is considered with the tensor product algebra structure. Also showing that this is coassociative, is straightforward. It follows from the formulas for $T$, given in Proposition 2.4 and Definition 2.7. This coproduct will be regular when also the assumption for $T^{\text{op}}$ is fulfilled. See [De2] for more details.
\nl
However, as we mentioned before, this is not (really) what we want. Eventually, we want to make the smash product algebra $A\# B$ into a multiplier Hopf algebra and so we need to define this coproduct on the smash product. Now, it takes again little effort to show that the result in the lemma will allow to define the linear map $\Delta$ from $A\# B$ to the multiplier algebra $M((A\# B) \ot (A\# B))$ (no matter what the right action of $A$ on $B$ is) and to show that it is coassociative. 
\snl
It is fair to say that this coproduct needs an algebra structure on $A\ot B$, but that it is essentially independent of the choice of this algebra structure.
\snl
It is in the next section that we will see what kind of compatibility conditions between the action and the coaction are needed for this coproduct to be an algebra homomorphism and not merely a coassociative linear map. Then, the approach will become different from what is obtained in [De2]. This will be explained.
\nl
Also the behaviour of this coproduct with respect to the $^*$-operation (in the case of multiplier Hopf $^*$-algebras), can not be treated here properly.
\snl
If we look at the condition on the action, necessary to have that the smash product $A\# B$ becomes a $^*$-algebra with the obvious $^*$-algebra structure (see a remark following Proposition 1.6 in the previous section), if we translate this to the case of a left $D$-module algebra $C$ and finally, if we dualize this condition, we arrive at the natural requirement that
$$\Gamma(S_A(a)^*)=((\iota_B\ot S_A)\Gamma(a))^*$$
for all $a\in A$.
\snl
This condition however is not sufficient to ensure that $\Delta$ is a $^*$-map on $A\# B$. We need one of the extra conditions on $\Gamma$ (namely the one that gives $\Gamma(aa')$ in terms of $\Gamma(a)$ and $\Gamma(a')$ for $a,a'\in A$ - see Section 3). 
\snl
The only thing we can obtain is that the map $\sigma T (S_A(\,\cdot\,)^* \ot S_B(\,\cdot\,)^*)$ is involutive from $A\ot B$ to itself (where $\sigma$ is the flip). This is no surprise because in the previous section, the condition relating the $^*$-structure with the action implies that $\sigma R ((\,\cdot\,)^* \ot (\,\cdot\,)^*)$ is involutive on $B\ot A$. And whereas this last property is sufficient to get that $A\# B$ is a $^*$-algebra, the first one is not sufficient to have that $\Delta$ is a $^*$-map. This is related with the fact that the antipode on $A\# B$ that we will get in Section 4, is not simply given by the antipode on $A$ and the antipode on $B$.
\snl
We come back to this problem in Section 3.  
\nl
\it Other cases and examples \rm
\nl
Now, we will do a few more things in this section, as we also did in the previous one. We will (very briefly) consider the case of a right $C$-comodule coalgebra $D$, we consider some obvious special cases  and we will look at the example.
\snl
Here is the definition of a right coaction. We will use the same symbols for the objects associated with this case as explained already earlier. 

\inspr{2.9} Definition \rm
Let $C$ be a regular multiplier Hopf algebra and $D$ a vector space. Now consider the vector space $M_0(D\ot C)$ of 'elements' $x$ such that $x(1\ot c)$ and $(1\ot c)x$ are in $D\ot C$ for all $c\in C$.
Assume that $\Gamma:D\to M_0(D\ot C)$ is an injective linear map so that $(\Gamma \ot \iota_C)\Gamma = (\iota_D\ot \Delta_C)\Gamma$. Then $\Gamma$ is called a {\it right coaction } of $C$ on $D$. \hfill$\square$
\einspr

We have the same remarks as about left coactions. Now, the concept is in duality with that of a right action of $A$ on $B$ 
(when $A$ is paired with $C$ and $B$ with $D$ as before). The relevant formula is 
$$\langle b\ot a, \Gamma(d)\rangle = \langle b\tl a, d\rangle$$
for all $a\in A$ and $b\in B$.
\snl
The other formulas can be obtained from the ones for a left coaction using the transformation tool described earlier. They can also be obtained by dualizing formulas about a right action.
\snl
The cotwist map $T$ associated with this coaction is the map from $C\ot D$ to $D\ot C$, defined by 
$T(c\ot d)=(1\ot c)\Gamma(d)$. In terms of this map, coassociativity is expressed as
$$(\iota_D\ot\Delta_C)T=(T\ot \iota_C)(\iota_C\ot T)(\Delta_C\ot\iota_D) \quad \text{on}\quad C\ot D.$$

\iinspr{2.10} Definition \rm 
If also $D$ is a regular multiplier Hopf algebra and if
$$(\Delta_D\ot \iota_C)T=(\iota_C\ot T)(T\ot \iota_D)(\iota_C\ot \Delta_D)$$
on $C\ot D$, then $D$ is called a {\it right $C$-comodule coalgebra}. \hfill$\square$
\einspr

Again, we have to assume the natural equivalent of the Assumption 2.6 for this last formula as discussed before.
\snl
The reader should be aware of the similarity of the two formulas for the cotwist map in the case of a left comodule coalgebra and that of a right comodule coalgebra. The order of the equations is different, just as in the case of left and right module algebras (see a remark in Section 1).
\snl
The formula for the coproduct is the same. In this case, it is defined formally on $C\ot D$  by
$$\Delta(c\ot d)= (\iota_C \ot T \ot\iota_D)(\Delta_C(c) \ot \Delta_D(d)).$$
The same remarks as for the case of a left $B$-comodule coalgebra $A$ hold here.
\nl
Before we come to the example, let us say a few words about the case of a trivial coaction. The left coaction $\Gamma$ of $B$ on $A$ is trivial if $\Gamma(a)=1_B\ot a$ for all $a\in A$. This implies that the cotwist map $T$ is simply the flip map and that the coproduct is nothing else but the tensor coproduct. A similar remark holds for a trivial right coaction $\Gamma$ of $C$ on $D$.

\iinspr{2.11} Example \rm Let $H$ and $K$ be groups.
\snl
i) Now first assume that $K$ acts on $H$ from the right (as a group on a set) and use $h\tl k$ to denote the action of the element $k\in K$ on the element $h\in H$. Assume that the action is unital. Let $A$ be the group algebra $\Bbb C H$ of $H$ and let $B$ be the function algebra $F(K)$ of complex functions with finite support on $K$. Define a map $\Gamma: A \to M(B\ot A)$ by
$$\Gamma(h)(\delta_k \ot 1)= \delta_k \ot (h\tl k)$$
where $h\in H$,  $k\in K$ and where $\delta_k$ denotes the function on $K$ that is $1$ in $k$ and $0$ everywhere else. One easily verifies that $\Gamma$ is a left coaction of $B$ on $A$ (in the sense of Definition 2.1) and that it makes $A$ into a left $B$-comodule coalgebra (as in Definition 2.7). It is in duality with the left $D$-module algebra $C$ as constructed in Example 1.8.ii (in the sense of Remark 2.3). 
\snl
The expression for the coproduct on $A\ot B$ is given by
$$\Delta(h\ot \delta_k)=\sum_{k'k''=k} h \ot \delta_{k'} \ot (h \tl k') \ot \delta_{k''}$$
when $h\in H$ and $k\in K$. The summation is taken over elements $k',k''$ in $K$.
\snl
ii) Next, assume that $H$ acts on $K$ from the left and use $h\tr k$ to denote the action. Again assume that the action is unital. Let $C$ be the function algebra $F(H)$ and  $D$ the group algebra $\Bbb C K$. Define $\Gamma: D \to M(D\ot C)$ by
$$(1\ot \delta_h)\Gamma(k)=(h\tr k) \ot \delta_h$$
for all $h\in H$ and $k\in K$. Then we have a right coaction of $C$ on $D$ (as in Definition 2.9), making $D$ into a right $C$-comodule coalgebra (as in Definition 2.10). It is in duality with the example in 1.8.i. 
\snl
Now the coproduct on $C\ot D$ is given by the expression
$$\Delta(\delta_h\ot k)=\sum_{h'h''=h} \delta_{h'} \ot (h''\tr k) \ot \delta_{h''} \ot k$$
for  $h\in H$, $k\in K$ and where again the summation is taken over elements $h',h''\in H$. 

\hfill$\square$
\einspr

We also have natural $^*$-algebras structures and they are compatible with the coactions. These two cases are much simpler than the general case because $A$ and $D$ are Hopf algebras. Finally, if the group actions are trivial, then also the coactions are trivial and the cotwist maps are simply the flip maps. The coproduct becomes the tensor coproduct. 

\newpage

\bf 3. The smash coproduct on the smash product \rm
\nl
In this section, we start with a pair of two regular multiplier Hopf algebras $A$ and $B$. We assume that $B$ is a right $A$-module algebra (as in Section 1) and that $A$ is a left $B$-comodule coalgebra (as in Section 2). We will freely use the notations of the two previous sections.
\snl
We consider the smash product $AB$ (as reviewed in Section 1) and we will consider the coproduct $\Delta$ on $AB$ as already discussed in Section 2. We have seen that the algebra structure on $AB$ is needed because the coproduct does not map $AB$ into $AB\ot AB$, but rather in $M(AB\ot AB)$. On the other hand, we also saw that this coproduct is not really dependent on the algebra structure. 
\snl 
In the previous section, we only considered this coproduct as a coassociative linear map. In this section, we will consider it as a coproduct on the smash product. We will see what kind of compatibility conditions are needed between the action and the coaction for this coproduct to be an algebra homomorphism. It then will take little effort to see that we get a regular multiplier Hopf algebra, but this will be done in the next section. The conditions are of course the natural generalizations of the well-known conditions imposed by Majid (see e.g.\ Theorems 6.2.2 and 6.2.3 in [M3]). Following the scope of this paper, we will spend some more time to motivate these conditions, discuss different forms of it and again concentrate on the problem of coverings. 
\nl
Before we start this in a systematic way, let us first rewrite the formula for the coproduct $\Delta$, as given in the previous section (see the two formulas before Lemma 2.8). Using $AB$ to denote the smash product, we can write formally
$$\Delta(ab)=\sum_{(a)}(a_{(1)}\ot 1)\Gamma(a_{(2)})\Delta(b)$$
whenever $a\in A$ and $b\in B$. Now, we will first consider this coproduct on $B$, then on $A$ and finally, we will consider the commutation rules so that, when defined on the algebra $AB$, we get a homomorphism. In the process of this way of constructing the coproduct, we will be concerned with the various conditions that are necessary and sufficient for this map to be an algebra map. We will, in order to avoid notational conflicts, use $\Delta_\#$ for the coproduct on $AB$ in stead of merely $\Delta$ as we did in the previous section.
\nl
\it The coproduct $\Delta_\#$ on the algebras $A$ and $B$ \rm
\nl 
We begin with the easiest part.

\inspr{3.1} Proposition \rm
There is a non-degenerate homomorphism $\Delta_\#: B \to M(AB\ot AB)$ given by 
$$\Delta_\#(b)=\Delta(b).$$
\hfill $\square$
\einspr

This is a triviality. We know e.g.\ that $(b'\ot 1)\Delta(b)$ and $\Delta(b)(1\ot b')$ are in $B\ot B$ for all $b,b'\in B$ and also because $AB=BA$, this easily implies that $\Delta_\#(b)\in M(AB\ot AB)$ for all $b$. Observe that the injection of $B\ot B$ in $AB\ot AB$ is a non-degenerate homomorphism so that $M(B\ot B)$ is sitting inside $M(AB\ot AB)$ (see e.g.\ the results in Proposition 1.5 in Section 1). And because  $\Delta$ is a non-degenerate homomorphism from $B$ to $M(B\ot B)$, it follows that $\Delta_\#$ will be a non-degenerate homomorphism from $B$ into $M(AB\ot AB)$.  
\snl
Coassociativity of $\Delta_\#$ on $B$ is an immediate consequence of coassociativity of $\Delta$ on $B$. We also get that 
$$\Delta_\#(B)(1\ot B)=B\ot B \qquad\qquad \text{and} \qquad\qquad (B\ot 1)\Delta_\#(B)=B\ot B$$
as sitting inside the multiplier algebra $M(AB\ot AB)$.

\nl
The next step is already less obvious. First, we define $\Delta_\#$ on $A$ and we prove some elementary properties.

\inspr{3.2} Proposition \rm
There is a linear map $\Delta_\#:A\to M(AB\ot AB)$ defined by 
$$\Delta_\#(a)=\sum_{(a)}(a_{(1)}\ot 1)\Gamma(a_{(2)}).$$
We have 
$$ (AB\ot 1)\Delta_\#(A)=AB\ot A \qquad \text{and} \qquad
         \Delta_\#(A)(B\ot A)=AB \ot A.$$

\bf \snl Proof: \rm
i) First we will show that $\Delta_\#(a)$, as in the formulation of the proposition, is a well-defined element in $M(AB\ot AB)$. We follow an argument that we gave already in the previous section. First multiply with $a'\ot 1$ from the left (with $a'\in A$). We get
$$\sum_{(a)}(a'a_{(1)}\ot 1)\Gamma(a_{(2)})\subseteq (A\ot 1)\Gamma(A)$$
and because $(B\ot 1)\Gamma(A)\subseteq B\ot A$, we see already that 
$$(AB\ot 1)\Delta_\#(a)\subseteq (AB\ot 1)\Gamma(A)\subseteq AB\ot A.$$
We have looked closer at this argument in  Example 3.5 of [VD4].
\snl  
Next, multiply with $b\ot a'$ from the right (with $a'\in A$ and $b\in B$).
We get 
$$\sum_{(a)}(a_{(1)}\ot 1)\Gamma(a_{(2)})(b\ot a')\subseteq AB\ot A$$
because of Assumption 2.6 (as discussed in the previous section). So we see that also 
$\Delta_\#(a)(b\ot a')\in AB\ot A$. The two results together will give that $\Delta_\#(a)$ is well-defined in $M(AB\ot AB)$. This proves the first part of the proposition.
\snl
ii) We will now show that this map satisfies the equalities. The arguments are close to the ones used in i). Indeed, if we look at the different steps in the first part of the proof above, we see that actually $(AB\ot 1)\Delta_\#(A)=AB\ot A$. We need  to use that $B\ot A=(B\ot 1)\Gamma(A)$ and that $A\ot A=(A\ot 1)\Delta(A)$. On the other hand, when we look at the other side, we need the equality
$$A\ot B\ot A=((\iota_A\ot \Gamma)\Delta_A(A))(1\ot B \ot A)$$
because then we have $AB\ot A=\Delta_\#(A)(B\ot A)$. Now, this equality is an easy consequence of the basic Assumption 2.6 and the one in Definition 2.7.  \hfill $\square$
\einspr

Whereas coassociativity of $\Delta_\#$ on $B$ was more or less obvious, this is not the case for $\Delta_\#$ on $A$. It is possible to argue (or verify) that $\Delta_\#$ is coassociative on $A$, but for this, we need $\Delta_\#$ on $AB$. This has been argued already in the previous section. We will come back to this problem later in this section (see Theorem 3.14) and in the next section (see a remark following the proof of Proposition 4.2).
\snl
We now look for necessary and sufficient conditions to ensure that $\Delta_\#(aa')=\Delta_\#(a)\Delta_\#(a')$ for all $a,a'\in A$. We will first prove the following result.

\inspr{3.3} Proposition \rm
 Define a linear map $P:B \ot A \to B\ot A$ by
$$P(b\ot a)= \sum_{(a)} ((b\tl a_{(1)})\ot 1)\Gamma(a_{(2)}).$$
Then $\Delta_\#(aa')=\Delta_\#(a)\Delta_\#(a')$ for all $a,a'\in A$  if and only if 
$$(\iota_B\ot m_A)P_{13}P_{12}=P(\iota_B\ot m_A)$$
on $B\ot A\ot A$. \hfill $\square$
\einspr

Before we prove this proposition, let us first make a few remarks about this map $P$ and the condition imposed on it.
\snl
First notice that the map $P$ can be written as a composition $T^{\text{op}}R^{\text{op}}$ where
$$\align R^{\text{op}}(b\ot a)&=\sum_{(a)} a_{(2)} \ot (b \tl a_{(1)}) \\
         T^{\text{op}}(a\ot b)&=(b\ot 1) \Gamma(a).
\endalign$$
Compare with the maps $T$ and $R$, introduced earlier, for understanding the notation. It follows from earlier arguments that $P$ is well defined. We also have a case of iterated coverings (cf.\ Example 3.5 in [VD4]).
\snl
Also compare the equation satisfied by $P$ with the first formula in Proposition 1.3. It is very similar but slightly different. The difference is mainly due to the fact that $P$ maps $B\ot A$ to itself whereas $R$ maps this space into $A\ot B$. We will also come back to this later (see a remark after Assumption 3.9 below).
\snl
For a better understanding of the result and the proof, as well as for results to come, it turns out to be useful to introduce the following.

\inspr{3.4} Notation \rm
Consider the right action of $AB$ on $B$, given before in Section 1 (see a remark following Proposition 1.6), by 
$$y\bullet (ab)=(y\tl a)b$$
whenever $a\in A$ and $b,y\in B$. Combine it with right multiplication of $A$ on itself to a right action of $AB\ot A$ on $B\ot A$. So
$$(y\ot x)\bullet (c\ot a)=(y\bullet c) \ot xa$$
when $a,x\in A$, $y\in B$ and $c\in AB$. \hfill $\square$
\einspr

This action is not necessarily faithful, but it is unital and so we can extend it to the multiplier algebra $M(AB\ot A)$. In particular, we can consider the action of the elements $\Delta_\#(a)$ for any $a\in A$. Then the following can be shown.

\inspr{3.5} Lemma \rm
We have $P(b\ot a) = (b\ot 1)\bullet \Delta_\#(a)$ for all $a\in A$ and $b\in B$.

\snl\bf Proof\rm: The formula has to be read correctly as
$$(1\ot x)P(b\ot a) = (b\ot x)\bullet \Delta_\#(a)$$
for all $x\in A$. As we clearly have 
$$(b\ot x)\bullet \Delta_\#(a)=\sum_{(a)}((b\tl a_{(1)})\ot x)\Gamma(a_{(2)})$$
whenever $a,x\in A$ and $b\in B$, the result is true. \hfill $\square$ 
\einspr

With this result, it is fairly easy to show one direction of Proposition 3.3. This is the content of the following lemma.

\inspr{3.6} Lemma \rm
If $\Delta_\#(aa')=\Delta_\#(a)\Delta_\#(a')$ for all $a,a'\in A$  then 
$$(\iota_B\ot m_A)P_{13}P_{12}=P(\iota_B\ot m_A)\qquad \text{on} \qquad B\ot A\ot A.$$ 

\snl\bf Proof: \rm
For all $a,a'\in A$ and $b\in B$ we have
$$\align (\iota_B\ot m_A)P_{13}P_{12}(b\ot a\ot a')
                        &=(\iota_B\ot m_A)P_{13}(((b\ot 1)\bullet \Delta_\#(a))\ot a')\\
                        &=((b\ot 1)\bullet\Delta_\#(a))\bullet\Delta_\#(a')\\
                        &=(b\ot 1)\bullet (\Delta_\#(a)\Delta_\#(a'))
\endalign$$  
and clearly also 
$$P(b\ot aa')=(b\ot 1)\bullet \Delta_\#(aa').$$ 
Strictly speaking, we should multiply these equations all the time with an element of $A$ from the left in the second factor. \hfill $\square$
\einspr

If on the other hand we have the condition on $P$, we see from the proof that 
$$(b\ot 1)\bullet \Delta_\#(aa')=(b\ot 1)\bullet (\Delta_\#(a)\Delta_\#(a'))$$
for all $b\in B$. We cannot conclude that $\Delta_\#(aa')=\Delta_\#(a)\Delta_\#(a')$ for all $a,a'\in A$
as the action of $AB\ot A$ on $B\ot A$ is not necessarily faithful. To prove the converse in Proposition 3.3, we have to use the action of $AB$ on $AB$, given by right multiplication. This is faithful as the product in $AB$ is non-degenerate. This will be done in the next lemma. 

\inspr{3.7} Lemma \rm
If $(\iota_B\ot m_A)P_{13}P_{12}=P(\iota_B\ot m_A)$ then $\Delta_\#(aa')=\Delta_\#(a)\Delta_\#(a')$ for all $a,a'\in A$.

\snl\bf Proof: \rm
Let $a,x\in A$ and $y\in B$. Then 
$$\align(xy\ot 1)\Delta_\#(a) &= \sum_{(a)}(xya_{(1)}\ot 1)\Gamma(a_{(2)})\\
                        &= \sum_{(a)}(xa_{(1)}\ot 1)((y\ot 1)\bullet\Delta_\#(a_{(2)})).
\endalign$$
If we replace in this equation $a$ by the product $aa'$, and if we assume the condition on $P$, it follows from the remark just made before this lemma that
$$\align (xy\ot 1)\Delta_\#(aa')
      &= \sum_{(a)(a')} (xa_{(1)}a'_{(1)}\ot 1)((y\ot 1)\bullet \Delta_\#(a_{(2)}a'_{(2)}))\\
      &= \sum_{(a)(a')} (xa_{(1)}a'_{(1)}\ot 1)((y\ot 1)\bullet\Delta_\#(a_{(2)})\bullet\Delta_\#(a'_{(2)}))\\
      &= \sum_{(a)} (xa_{(1)}\ot 1)((y\ot 1)\bullet\Delta_\#(a_{(2)}))\Delta_\#(a')\\
      &= (xy\ot 1)\Delta_\#(a)\Delta_\#(a').
\endalign$$
As the action of $AB$ on $AB$ is faithful, we get $\Delta_\#(aa')=\Delta_\#(a)\Delta_\#(a')$. \hfill $\square$
\einspr

Combining Lemma 3.6 with Lemma 3.7 we obtain Proposition 3.3. 

\snl
Before we continue, let us look at another form of the condition that looks more like the one used in Hopf algebra theory. We can also prove the following result but the precise  coverings are quite involved. Moreover, the result is not really important for our approach. Nevertheless, let us consider it briefly. It is a good illustration of how complicated these coverings can be.

\inspr{3.8} Proposition \rm
We have 
$\Delta_\#(aa')=\Delta_\#(a)\Delta_\#(a')$ for all $a,a'\in A$ 
if and only if
$$\Gamma(aa')=\sum_{(a)(a')}((a_{(-1)}\tl a'_{(1)})\ot a_{(0)})\Gamma(a'_{(2)})$$
for all $a,a'\in A$.

\snl\bf Proof: \rm
Of course, one first has to give a meaning to the second formula by using the right coverings. There are several possibilities. The easiest one is obtained if we multiply with an element of $B$ in the first factor from the left.
\snl
Then the left hand side is alright. For the right hand side, we have a part of the formula like
$$b(a_{(-1)}\tl a'_{(1)})=((b\tl S(a'_{(1)}))a_{(-1)})\tl a'_{(2)}$$
and we see that first $b$ will cover $a'_{(1)}$ (through the action), then $a_{(-1)}$ will be covered (by multiplication) and finally, $a'_{(2)}$ will again be covered (through the action). So, a modified form of the right hand side of this formula for $\Gamma(aa')$ is covered (by multiplication with an element of $B$ from the left in the first factor). Another way to fully cover the expression is by also multiplying from the right with elements of the form $b'\ot a''$ where $a''\in A$ and $b'\in B$. Then Assumption 2.6 is used, but still, things get quite complicated.
\snl
To prove the result, we actually need to use this last covering together with the following arguments.
\snl 
The formula can be rewritten as $\Gamma(aa')=\Gamma(a)\bullet \Delta_\#(a')$. This is correct for all $a,a'$ if and only if 
$$\sum_{(a)(a')}(a_{(1)}a'_{(1)}\ot 1)\Gamma(a_{(2)}a'_{(2)})
=\sum_{(a)(a')}(a_{(1)}a'_{(1)}\ot 1)(\Gamma(a_{(2)})\bullet \Delta_\#(a'_{(2)}))$$
for all $a,a'\in A$. But the right hand side of this equation is equal to 
$$\sum_{(a)}(a_{(1)}\ot 1)\Gamma(a_{(2)}) \Delta_\#(a')=\Delta_\#(a)\Delta_\#(a')$$
and the result follows. \hfill $\square$
\einspr

We see that things get quite complicated if we want to do this completely rigorously? 
In any case, when we would have Hopf algebras, no  
problem would occur and the proof of the result in Proposition 3.8 becomes more or less obvious.
  
\nl
\it The coproduct $\Delta_\#$ on $AB$ \rm
\nl
Next, we want to have that $\Delta_\#$ is a homomorphism on the smash product $AB$. Since we have already shown that it is a homomorphism on $A$ and on $B$ (provided we assume the right conditions), it remains to verify that the commutation rules 
$ba=\sum_{(a)}a_{(1)}(b\tl a_{(2)})$ are respected. In other words, we need to have
$$\Delta_\#(b)\Delta_\#(a)=\sum_{(a)}\Delta_\#(a_{(1)})\Delta_\#(b\tl a_{(2)})$$
in $M(AB \ot AB)$ for all $a\in A$ and $b\in B$.
This requires 3 steps. We will need a way to commute elements of $B$ with the first factor of $\Delta_\#(a)$, a way to commute elements of $B$ with the second factor of $\Delta_\#(a)$ and finally, we will need a formula for $\Delta(b\tl a)$.
\snl
We begin with the last property. Indeed, we have seen before how the coaction $\Gamma$ behaves with respect to the product in $A$. What we need here is how the coproduct of $B$ relates with the action. This is precisely the dual situation. We know how to deduce one from the other.
\snl
We start with the formula, given in Proposition 3.3. We had 
$$(\iota_B\ot m_A)P_{13}P_{12}=P(\iota_B\ot m_A)$$
on $B\ot A\ot A$ where $P$ is defined from $B\ot A$ to itself by
$$P(b\ot a)= \sum_{(a)} ((b\tl a_{(1)})\ot 1)\Gamma(a_{(2)}).$$
Next, we transform these formulas to the case of a left action of $D$ on $C$ and a right coaction $\Gamma$ of $C$ on $D$ (cf.\ Definitions 2.9 and 2.10). Then we get the map $P'$ from $D\ot C$ to itself, given by
$$P'(d\ot c)=\sum_{(d)}\Gamma(d_{(1)})(1\ot (d_{(2)} \tr c)).$$ 
It should satisfy the equation
$$(m_D\ot \iota_C)P'_{13}P'_{23}=P'(m_D\ot \iota_C)$$
on $D\ot D\ot C$.
Finally, if we look for the adjoint of this map $P'$, we find
$$\align \langle b\ot a, P'(d\ot c) \rangle 
                &= \sum_{(d)} \langle b\ot a, \Gamma(d_{(1)})(1\ot (d_{(2)}\tr c))\rangle \\
                &= \sum_{(a),(d)} \langle (b\tl a_{(1)})\ot a_{(2)}, d_{(1)}\ot (d_{(2)}\tr c)\rangle \\
                &= \langle P(b\ot a),d\ot c \rangle
\endalign$$
for all $a,b,c$ and $d$ in $A,B,C$ and $D$ respectively. We see that $P'$ and $P$ are each others adjoints in this situation. Now, if we take the adjoint of the condition on $P'$, we arrive at the following natural assumption about $P$.

\inspr{3.9} Assumption \rm
$P_{23}P_{13}(\Delta_B\ot \iota_A)=(\Delta_B\ot \iota_A)P$ on $B\ot A$.
\hfill $\square$
\einspr

This formula should be compared with the formulas in Proposition 2.4 and Definition 2.7. The covering problem can be solved in two ways. Either we write the assumption as $P_{13}(\Delta_B\ot \iota_A)=P_{23}^{-1}(\Delta_B\ot \iota_A)P$ on $B\ot A$.
The left hand side is now covered if we multiply with elements of $B$ in the second factor (left or right). The right hand side is covered with elements of $B$ in the first factor (again left or right). This is similar as for the equation in Definition 2.7.
\snl
The other possibility goes as follows. Given $a\in A$ and $b\in B$, we see that the map  $q \mapsto (b\ot 1)P(q\ot a)$ from $B$ to $B\ot A$ has the variable covered from the left. Indeed,
$$(b\ot 1)P(q\ot a)=\sum_{(a)}(b(q\tl a_{(1)})\ot 1)\Gamma(a_{(2)})$$
and we have seen before  that $b$ and $a$ will eventually cover $q$ from the left. It follows that the equation in Assumption 3.9 is also well-covered if we multiply with an element of $B$ from the left, either in the first or in the second factor. 
\snl
Now, we are interested in another form of this equation. For this purpose, we let $AB$ act from the right on $B$ as we did before, but now we also consider $AB\ot AB$ as acting on $B\ot B$. We get the following:

\iinspr{3.10} Proposition \rm The assumption in 3.9 is fulfilled if and only if
 $\Delta(b\tl a)=\Delta(b)\bullet \Delta_\#(a)$ for all $a\in A$ and $b\in B$. This equation will be well-covered if we multiply from the left in the first factor with any element of $B$.

\snl\bf Proof: \rm
We start with the equation 
$$(\Delta_B\ot\iota_A )P(b\ot a)=P_{23}P_{13}(\Delta_B(b)\ot a)$$
with $a\in A$ and $b\in B$. We apply $\iota_B\ot \iota_B\ot \varepsilon_A$. On the left hand side, we get $\Delta_B(b\tl a)$ because $(\iota\ot \varepsilon_A)\Gamma(a)=\varepsilon_A(a)1$. For the right hand side we get
$$\align \sum_{(b)}(\iota_B\ot\iota_B\ot \varepsilon_A)&P_{23}P_{13}(b_{(1)}\ot b_{(2)}\ot a)\\
    &= \sum_{(a),(b)}(\iota_B\ot\iota_B\ot \varepsilon_A)P_{23}((b_{(1)}\tl a_{(1)})\ot b_{(2)}\ot 1)\Gamma_{13}(a_{(2)})\\
      &= \Delta(b)\bullet \Delta_\#(a)
\endalign$$
(where we again have used $(\iota_B\ot \varepsilon_A)P(b'\ot a')=b'\tl a'$). Remark that everything is well covered if we multiply with an element of $B$ from the left, either in the first or in the second factor.
\snl
Also conversely, if $\Delta(b\tl a)=\Delta(b)\bullet \Delta_\#(a)$ for all $a,b$, the assumption in 3.9 will be satisfied.
\hfill $\square$
\einspr
Now, we have taken care of the third step as mentioned before. Another step in the whole procedure is to commute elements of $B$ in the second factor with elements $\Delta_\#(a)$. In other words, we look at a formula for $(1\ot b)\Delta_\#(a)$ and we try to move $b$ to the other side. Because this amounts to commuting $B$ with $A$, there is no need for an extra condition here. This is what we get:

\iinspr{3.11} Proposition \rm
For all $a\in A$ and $b\in B$ we have
$$(1\ot b)\Delta_\#(a)=\sum_{(a)}\Delta_\#(a_{(1)})((1\ot b)\bullet \Gamma(a_{(2)})).$$
Here, we use $\bullet$ to denote the right action of $B\ot A$ on $B\ot B$, obtained from right multiplication by $B$ in the first factor and the right action of $A$ in the second factor. Again the action is extended to the multiplier algebras.

\snl \bf Proof: \rm
We first present a formal argument and we discuss the details about the necessary coverings later.
\snl
We know that 
$(\iota\ot\Delta_A)\Gamma(a)=\sum_{(a)}\Gamma_{12}(a_{(1)})\Gamma_{13}(a_{(2)})$.
Therefore 
$$\align (1\ot b)\Gamma(a)&=\sum_{(a)} a_{(-1)}\ot ba_{(0)}\\
                &=\sum_{(a)}a_{(-1)} \ot a_{(0)(1)}(b\tl a_{(0)(2)}) \\
                &=\sum_{(a)}a_{(1)(-1)}a_{(2)(-1)} \ot a_{(1)(0)}(b\tl a_{(2)(0)})\\
                &=\Gamma(a_{(1)})((1\ot b)\bullet\Gamma(a_{(2)})).
\endalign$$
If we replace $a$ by $a_{(2)}$ and multiply with $a_{(1)}$ from the left in the first factor, we get 
$$\sum_{(a)}(a_{(1)}\ot b)\Gamma(a_{(2)})=\sum_{(a)}(a_{(1)}\ot 1)\Gamma(a_{(2)})((1\ot b)\bullet \Gamma(a_{(3)}))$$
and this is the the required formula.
\snl
In the last formula, there is no problem with covering $a_{(1)}$. We simply multiply with any element of $A$ in the first factor from the left. This takes care of the last step in the argument above. 
\snl
Next, we multiply all expressions in the first series of formulas with elements of the form $b'\ot a'$, from the right, where $a'\in A$ and $b'\in B$. Then we know from earlier discussions that at all stages in the calculation, everything is well-defined in $B\ot A$ and that all equations hold.
\hfill $\square$ 
\einspr
Finally, we  need to take $b$ to the other side in the formula $(b\ot 1)\Delta_\#(a)$. This is fundamentally different from the case before with $b$ in the second factor because the first factor of $\Delta_\#(a)$ is in $AB$. Indeed, we need an extra assumption for this:

\iinspr{3.12} Assumption \rm
We assume that $T\circ R=T^{\text{op}}\circ R^{\text{op}}$.
\hfill $\square$
\einspr

Recall the formulas for the maps $R$ and $T$ (see Proposition 1.2 and Proposition 2.4):
$$\align R(b\ot a)&=\sum_{(a)} a_{(1)} \ot (b\tl a_{(2)}) \\
         T(a\ot b)&=\Gamma(a)(b\ot 1)
\endalign$$
where $a\in A$ and $b\in B$. Similarly $R^{\text{op}}$ and $T^{\text{op}}$ are defined:
$$\align R^{\text{op}}(b\ot a)&=\sum_{(a)} a_{(2)} \ot (b\tl a_{(1)}) \\
         T^{\text{op}}(a\ot b)&=(b\ot 1)\Gamma(a)
\endalign$$
with $a\in A$ and $b\in B$. The composition $T^{\text{op}}\circ R^{\text{op}}$ is nothing else but the operator $P$ as introduced in Proposition 3.3 while the composition of $T$ with $R$ is given by 
$$  (T\circ R)(b\ot a)=\sum_{(a)} \Gamma(a_{(1)})((b\tl a_{(2)})\ot 1).$$
Then we have the following. It is essentially a reformulation of this condition and we see from it that this condition indeed allows to move elements $b$ in the expression $(b\ot 1)\Delta_\#(a)$ to the other side.

\iinspr{3.13} Lemma \rm
For all $a\in A$ and $b\in B$ we have $$(b\ot 1)\Delta_\#(a)=\sum_{(a)}\Delta_\#(a_{(1)})((b\tl a_{(2)}) \ot 1).$$

\snl\bf Proof: \rm 
For all $a\in A$ and $b\in B$ we get 
$$\align  (b\ot 1)\Delta_\#(a)&= \sum_{(a)}(ba_{(1)} \ot 1) \Gamma(a_{(2)}) \\
&= \sum_{(a)} (a_{(1)}(b\tl a_{(2)}) \ot 1)\Gamma(a_{(3)})\\
&= \sum_{(a)} (a_{(1)} \ot 1) P(b\ot a_{(2)}).
\endalign$$
In this case, we can multiply with elements of $A$ in the first factor from the left to cover every expression properly.
\snl
Now we use the assumption and we get
$$\align (b\ot 1)\Delta_\#(a)
&= \sum_{(a)} (a_{(1)} \ot 1)\Gamma(a_{(2)})((b\tl a_{(3)})\ot 1) \\
&= \sum_{(a)}\Delta_\#(a_{(1)})(b\tl a_{(2)} \ot 1).
\endalign$$
Again, we multiply with an element of $A$ on the left in the first factor and everything will be covered. Observe that $b$ will cover $a_{(2)}$ and $a_{(3)}$ respectively through the action.
\hfill $\square$
\einspr

Now, we are ready to conclude with the main result of this section.

\iinspr{3.14} Theorem \rm
Assume that we have two regular multiplier Hopf algebras $A$ and $B$ and that $B$ is a right $A$-module algebra and that $A$ is a left $B$-comodule coalgebra as before. Consider the associated maps and assume that we have
$$\align P(\iota_B\ot m_A)&= (\iota_B\ot m_A)P_{13}P_{12} \quad \text{on} \quad B\ot A\ot A \\
       (\Delta_B\ot \iota_A)P&= P_{23}P_{13}(\Delta_B\ot \iota_A) \quad \text{on} \quad B\ot A \\
       \endalign$$ 
as well as $T\circ R =T^{\text{op}}\circ R^{\text{op}}$  (i.e.\ $P=T\circ R $) on $B\ot A$.
\snl
Then the coproduct $\Delta_\#$ is a (well-defined) homomorphism on the smash product $AB$.

\snl \bf Proof: \rm
We have defined $\Delta_\#$ on $A$ and on $B$ and we know that $\Delta_\#(ab)$ is equal to $\Delta_\#(a)\Delta_\#(b)$. We know from Section 2 that this map is coassociative. From the results in Proposition 3.1 and 3.2, it follows that it is non-degenerate.
Because of Proposition 3.3, we know that the first condition implies that $\Delta_\#$ is a homomorphism on $A$. And if we combine the results obtained using the other two conditions, we find for all $a\in A$ and $b\in B$ that
$$\align \Delta_\#(b)\Delta_\#(a) &= \sum_{(b)}(b_{(1)}\ot b_{(2)})\Delta_\#(a) \\
                &=\sum_{(a),(b)}(b_{(1)}\ot 1)\Delta_\#(a_{(1)}) ((1\ot b_{(2)})\bullet \Gamma(a_{(2)}))\\
                &=\sum_{(a)}\Delta_\#(a_{(1)})(\Delta(b)\bullet \Delta_\#(a_{(2)}))\\
                &=\sum_{(a)}\Delta_\#(a_{(1)})\Delta(b\tl a_{(2)}).
\endalign$$
We can cover properly if we multiply with an element of $AB$ in the first factor.
\snl
This proves that $\Delta_\#$ is a homomorphism on the smash coproduct $AB$. \hfill $\square$
\einspr
The reader should compare this result with Theorem 6.2.3 of [M3]. Remark that it follows from the second condition that $\varepsilon_B(b\tl a)=\varepsilon_B(b)\varepsilon_A(a)$ and so, there is no need to add this condition as an assumption. When dealing with Hopf algebras, the first condition will imply that $\Gamma(1_A)=1_A\ot 1_B$ (because $\Delta_\#(1)=1$). In our setting however, we still can say that $\Delta_\#(1)=1$ as it makes sense to extend $\Delta_\#$ to the multiplier algebra. However, it is not obvious how to extend $\Gamma$ itself to the multiplier algebra $M(A)$. See also a remark in the proof of Theorem 4.3 and in Section 5. Nevertheless, this explains why we do not need to impose an extra condition of this type on $\Gamma$ as is done in [M3].
\nl
\it Other cases and examples \rm
\nl
At this stage, we can complete the results considered in the previous sections about the $^*$-algebra case. We get the following.

\iinspr{3.15} Theorem \rm
Let $A$ and $B$ be multiplier Hopf $^*$-algebras. Assume that $B$ is a right $A$-module algebra as before and also that 
$$(b\tl a)^*=b^*\tl S(a)^*$$
for all $a\in A$ and $b\in B$. Then $AB$ is a $^*$-algebra for the involution defined by $(ab)^*=b^*a^*$. Assume also as before that $A$ is a left $B$-comodule coalgebra and that 
$$\Gamma(S_A(a)^*)=((\iota_B\ot S_A)\Gamma(a))^*$$
for all $a\in A$. If furthermore, we have the assumptions as in Theorem 3.14, then $\Delta_\#$ is a $^*$-homomorphism on the $^*$-algebra $AB$.

\snl\bf Proof: \rm
We have seen already in Section 1 that $AB$ becomes a $^*$-algebra with the obvious involution. We also have seen in Section 2 that the condition on the coaction is a natural one, but that it is not sufficient to ensure that $\Delta_\#$ is a $^*$-map. We will now see however that this condition, together with the property that $\Delta_\#$ is a homomorphism, will imply that it is a $^*$-homomorphism.
\snl
In order to show that this is the case, we just have to verify that we have a $^*$-map on the components $A$ and $B$ of $AB$. There is obviously no problem on $B$ because $\Delta_\#$ coincides with $\Delta_B$ on $B$ and by assumption, $\Delta_B$ is a $^*$-homomorphism. So, we just have to verify that $\Delta_\#(a^*)=\Delta_\#(a)^*$ for all $a\in A$. 
\snl
Using the formula for $\Gamma(S(a)^*)$ we get
$$\align \Delta_\#(S(a)^*) &= \sum_{(a)} (S(a_{(2)})^* \ot 1)\Gamma(S(a_{(1)})^*) \\
                      &= \sum_{(a)} ((\iota\ot S)\Gamma(a_{(1)})(S(a_{(2)}) \ot 1))^*
\endalign$$
for all $a$. So, if we want to show that $\Delta_\#(S(a)^*)=\Delta_\#(S(a))^*$ for all $a\in A$, it will be sufficient to prove that 
$$\sum_{(a)}(S(a_{(2)})\ot 1)\Gamma(S(a_{(1)}))=\sum_{(a)}(\iota\ot S)\Gamma(a_{(1)})(S(a_{(2)})\ot 1)$$
for all $a$. Observe that this equation does not involve the involutions anymore.
\snl
To prove this result, it is allowed to multiply with  $\Delta_\#(a_{(3)})$ from the right and prove the resulting equation. For the left hand side, we get
$$\sum_{(a)}\Delta_\#(S(a_{(1)})) \Delta_\#(a_{(2)})=\varepsilon(a)\Delta_\#(1)=\varepsilon(a)1$$
while for the right hand side, we obtain
$$\sum_{(a)}((\iota\ot S)\Gamma(a_{(1)}))\Gamma(a_{(2)}).$$
And indeed, these two expressions coincide for, if we start with the formula
$$(\iota\ot\Delta)\Gamma(a)=\sum_{(a)}\Gamma_{12}(a_{(1)})\Gamma_{13}(a_{(2)}),$$
apply $m(S\ot\iota)$ on the last two factors and use that $(\iota\ot\varepsilon)\Gamma(a)=\varepsilon(a)1$, we get the desired equality. We can cover by multiplying from the right with an element of $B$ in the first factor and an element of $A$ in the last factor. \hfill $\square$
\einspr

Observe that the extra condition which is needed here to show that $\Delta_\#$ is a $^*$-map is just that it is a algebra homomorphism on $A$. This is the condition in Proposition 3.3.
\nl
As we have done in Section 1 and Section 2, also here we will briefly consider the dual case. 
\snl
So, we have again two regular multiplier Hopf algebras $C$ and $D$. We assume that $C$ is a left $D$-module algebra (as in Section 1) and that $D$ is a right $C$-comodule coalgebra (as in Section 2). As before, we will also use $\Gamma$ to denote this coaction. In this case, the twist maps $R,R^{\text{op}}:D\ot C \to C \ot D$ and the cotwist maps $T,T^{\text{op}}:C\ot D \to D\ot C$ are given by
$$\align R(d\ot c) &= \sum_{(d)} (d_{(1)}\tr c) \ot d_{(2)} \\
         R^{\text{op}}(d\ot c) &=\sum_{(d)} (d_{(2)}\tr c) \ot d_{(1)} \\
         T(c\ot d) &= (1\ot c)\Gamma(d) \\
         T^{\text{op}}(c\ot d) &= \Gamma(d)(1\ot c).
\endalign$$
Similarly, the map $P:D\ot C \to D\ot C$, defined as $T^{\text{op}}\circ R^{\text{op}}$, in this case is given by
$$P(d\ot c) = \sum_{(d)} \Gamma(d_{(1)})(1\ot (d_{(2)}\tr c)).$$
\snl
The dual version of Theorem 3.14 is now the following:

\iinspr{3.16} Theorem \rm
Let $C$ and $D$ be two regular multiplier Hopf algebras. Assume that $C$ is a left $D$-module algebra and that $D$ is a right $C$-comodule coalgebra. With the notations as above, assume furthermore that 
$$\align 
    P(m_D\ot \iota_C) &= (m_D\ot \iota_C)P_{13}P_{23} \quad \text{on} \quad D\ot D\ot C \\
   (\iota_D\ot \Delta_C) P &= P_{12}P_{13}(\iota_D\ot\Delta_C) \quad \text{on} \quad D\ot C
\endalign$$
as well as $P=T\circ R$ on $D\ot C$.
\snl
Then the smash coproduct $\Delta_\#$, defined as  
$\Delta_\#(cd)=\sum_{(d)}\Delta_C(c)\Gamma(d_{(1)})(1\ot d_{(2)})$
is an algebra homomorphism on the smash product $CD$.
\hfill $\square$
\einspr

It is also clear from all the considerations above that, when $A$ is paired with $C$ and when $B$ is paired with $D$ and if actions and coactions are adjoint to each other, then the tensor product pairing between $AB$ and $CD$ will give a pairing in the sense that the product on one algebra is adjoint to the coproduct on the other one. In the next section, we will see that $AB$ and $CD$ are regular multiplier Hopf algebras and that we do have a pairing in the sense of multiplier Hopf algebras.
\snl
Remark that the first condition  in Theorem 3.16 is like the first condition in Theorem 3.14 that we get when we apply the standard rules to pass from the pair $(A,B)$ to the pair $(C,D)$. Similarly for the other conditions. However, in the dual pair picture, the first and the second conditions in Theorem 3.16 are dual to the second and the first conditions  in Theorem 3.14 respectively. The last condition is a 'self-dual' condition. It means that in the case of such a pairing, it is sufficient to impose one set of condition, the other set will follow from duality.
\snl
Also in this case, we can see what happens when we have multiplier Hopf $^*$-algebras. The relevant conditions relating the action and coaction with the involution become 
$$(d\tr c)^*=S(d)^*\tr c^* \qquad\text{ and }\qquad \Gamma(S(d)^*)=((S\ot \iota)\Gamma(d))^*$$
for all $c\in C$ and $d\in D$. Then, we get that the smash product is a $^*$-algebra and that the coproduct is a $^*$-homomorphism.
\nl
Before we consider the basic example, let us have a look at some {\it special cases} of Theorem 3.14.
\snl
First consider the {\it case of a trivial action} of $A$ on $B$. Then we see from Proposition 3.8 that $\Delta_\#$ will be an algebra map on $A$ if and only if $\Gamma(aa')=\Gamma(a)\Gamma(a')$ for all $a,a'\in A$. This means that $\Gamma$ is an algebra map and we have that $A$ is a left $B$-comodule {\it bi}-algebra. This then takes care of the first condition in Theorem 3.14. The second condition turns out to be nothing else but the comodule property $(\Delta_B \ot \iota_A)\Gamma=(\iota_B\ot \Gamma)\Gamma$. Finally, the last condition becomes $T=T^{\text{op}}$ as both $R$ and $R^{\text{op}}$ are the flip map. So, we need
$$\Gamma(a)(b\ot 1)=(b\ot 1)\Gamma(a)$$
for all $a\in A$ and $b\in B$. This is indeed equivalent with the property that $\Delta_\#(a)$ and $\Delta_\#(b)$ will commute for all $a,b$ (given that the action is trivial so that $A$ and $B$ commute in the smash product).
\snl
Of course, $\Gamma(a)(b\ot 1)=(b\ot 1)\Gamma(a)$ for all $a,b$ will be fulfilled if $B$ is abelian. Then we arrive at the situation of Theorem 1.15 in [De2]. Also if the coaction is trivial, this will be true. And it is not hard to imagine that there are also cases in between these two extremes where the condition will be satisfied.
\snl
Next, consider the case where the action is not necessarily trivial, but where we have a {\it trivial left coaction} of $B$ on $A$. Then we not only have that $\Delta_\#=\Delta_B$ on $B$ but also $\Delta_\#=\Delta_A$ on $A$. So, $\Delta_\#$ will automatically be an algebra map on $A$ also. This takes care of the first condition in Theorem 3.14. The second condition reads now
$$\Delta_B(b\tl a)=\sum_{(a)(b)}(b_{(1)}\tl a_{(1)}) \ot (b_{(2)}\tl a_{(2)})$$
for all $a\in A$ and $b\in B$. This means that $B$ is a right $A$-module {\it bi}-algebra. The third condition says that $R=R^{\text{op}}$ as $T$ and $T^{\text{op}}$ are the flip map. So, we  need
$$\sum_{(a)} a_{(1)} \ot (b\tl a_{(2)})=\sum_{(a)} a_{(2)} \ot (b\tl a_{(1)})$$
for all $a$ and $b$. This will be true if $A$ is cocommutative and then we arrive at the situation of Theorem 1.6 in [De1]. Of course, if also the action is trivial, the condition is again satisfied, but there may be other cases between these two extremes where this is also true.
\nl
Finally, we consider the example of a matched pair of groups.

\iinspr{3.17} Examples \rm
i) Consider a group $G$ with two subgroups $H$ and $K$ so that $G=KH$ and $H\cap K=\{e\}$. As explained already in the introduction, we get a left action $\tr$ of the group $H$ on the set $K$ and a right action $\tl$ of the group $K$ on the set $H$ defined by the formula
$$hk=(h\tr k)(h\tl k)$$
for $h\in H$ and $k\in K$.
\snl
As before, we use $A$ to denote the group algebra $\Bbb C H$ and $D$ for the group algebra $\Bbb C K$. And we use $C$ for the algebra $F(H)$ of functions with finite support on $H$ and $B$ for the algebra $F(K)$ of functions with finite support on $K$. The left action of $H$ on $K$ induces a right action of $A$ on $B$ making $B$ into a right $A$-module algebra. Similarly, the right action of $K$ on $H$ induces a left action of $D$ on $C$ making $C$ into a left $D$-module algebra. See Example 1.8 for details.
\snl
Also the right action of $K$ on $H$ makes $A$ into a left $B$-comodule coalgebra whereas the left action of $H$ on $K$ makes $D$ into a right $C$-comodule coalgebra. See Example 2.11 for details.
\snl
ii) The product on the smash product $AB$ is given by the formula
$$(hf)(h'f')=(hh')((f\tl h')f')$$
and the coproduct $\Delta_\#$ on $AB$ is given by
$$\Delta_\#(h\delta_k)=\sum_{k'k''=k} h\delta_{k'} \ot (h\tl k')\delta_{k''}$$
where $h,h'\in H$, $k,k',k''\in K$ and $f,f'\in F(K)$ and where $\delta_k$ denotes the function on $K$ that is one on the element $k\in K$ and $0$ on all other elements. Similarly, the product on the smash product $CD$ is defined by
$$(fk)(f'k')=(f(k\tr f'))(kk')$$
and the coproduct $\Delta_\#$ on $CD$ is given by 
$$\Delta_\#(\delta_h k)=\sum_{h'h''=h} \delta_{h'}(h''\tr k) \ot \delta_{h''}k$$
where $h,h',h''\in H$, $k,k'\in K$ and $f,f'\in F(H)$. See again Example 1.8 and Example 2.11 for details.
\snl
iii) We will now verify the conditions in Theorem 3.14 and Theorem 3.16.
\snl
Of course, we can directly verify that $\Delta_\#$ is a homomorphism on the smash product $AB$ (as is done in e.g.\ [VD-W]). Similarly for $\Delta_\#$ on $CD$. This second case will also follow, either by symmetry (applying the rules to convert the pair $(A,B)$ to the pair $(C,D)$) or by duality. But as we treat this example here to illustrate the general theory, let us rather look at the three conditions in each of the theorems.
\snl
The easiest one is the third condition $T\circ R =T^{\text{op}}\circ R^{\text{op}}$. Because $A$ is cocommutative, we have $R=R^{\text{op}}$ and because $B$ is abelian, we have $T=T^{\text{op}}$. So, this condition is valid in Theorem 3.14. Similarly, or by duality, the condition is satisfied in Theorem 3.16. 
\snl
Next, consider the first condition of Theorem 3.14. For this, we have to verify that $\Delta_\#(h)\Delta_\#(h_1)=\Delta_\#(hh_1)$ for all $h,h_1\in H$. Now we have
$$\Delta_\#(h)=\sum_{k,k'} h\delta_{k'} \ot (h\tl k')\delta_{k''}$$
and similarly for $\Delta_\#(h_1)$. For the product we get
$$\align \Delta_\#(h)\Delta_\#(h_1) 
      &=\sum h \delta_{k'}h_1\delta_{k'_1} \ot (h\tl k')\delta_{k''}(h_1 \tl k'_1)\delta_{k''_1} \\
      &=\sum (hh_1)(\delta_{k'}\tl h_1)\delta_{k'_1} \ot (h\tl k')(h_1 \tl k'_1)(\delta_{k''}\tl (h_1\tl k'_1))\delta_{k''_1}
\endalign$$
where the summation is taken over all $k',k'',k'_1,k''_1\in K$. This should be equal to
$$\Delta_\#(hh_1)=\sum (hh_1) \delta_{k'_1} \ot ((hh_1)\tl k'_1)\delta_{k''_1}$$
with summation over all $k'_1,k''_1\in K$. Because  we have $(\delta_{k'}\tl h_1)\delta_{k'_1}$ non zero only if $k'=h_1\tr k'_1$ we find that these two expressions will be the same if and only if $(hh_1)\tl k'_1=(h\tl k')(h_1\tl k'_1)$ whenever  
$k'=h_1\tr k'_1$. This means that we need
$$(hh_1)k'_1=(h\tl(h_1\tr k'_1))(h_1\tl k'_1)$$
for all $k'_1$. This is obtained when we express the equality $(hh')k=h(h'k)$ (with $h'=h_1$ and $k=k_1$) and use the definition of the actions and the uniqueness of the decomposition. Indeed
$$\align (hh')k &= ((hh')\tr k)(hh')\tl k) \\
         h(h'k) &= h (h'\tr k)(h'\tl k) \\
                &= (h\tr (h'\tr k))(h \tl (h'\tr k))(h'\tl k)
\endalign$$
so that not only $(hh')\tr k=h\tr (h'\tr k)$ but also $(hh')\tl k=(h \tl (h'\tr k))(h'\tl k)$ for all $h,h'\in H$ and $k\in K$. This completes the argument and we have shown that $\Delta_\#$ is indeed a homomorphism on $A$, i.e.\ the first condition in Theorem 3.14 is satisfied.
\snl
In a similar way, we get the first condition of Theorem 3.16 from the equality $h(kk')=(hk)k'$. By duality, we find that this property will give the second condition in Theorem 3.14 and it is also possible to check this. Similarly, by duality, we get the second condition in Theorem 3.16 from the first one in 3.14.
\snl
So, we see that all conditions are fulfilled and that the coproducts $\Delta_\#$ are homomorphisms on both $AB$ and $CD$ for this example.
\hfill $\square$
\einspr

In the next section, we will complete these two basic examples and show that we indeed get regular multiplier Hopf ($^*$-)algebras as expected.
 
\nl\nl

\bf 4. The bicrossproduct for regular multiplier Hopf algebras \rm
\nl
In this section, we will formulate and prove the {\it main results}. Remark however that the most important work has been done already in the previous section.
\snl
As before, we have two regular multiplier Hopf algebras $A$ and $B$. The algebra $B$ is a right $A$-module algebra and $A$ is a left $B$-comodule coalgebra. We consider the smash product $AB$ and the coproduct $\Delta_\#$ on $AB$ as studied in the previous section. In particular, we assume that the relations between the $A$-module structure and the $B$-comodule structure, as formulated in Theorem 3.14, are fulfilled.
\snl
Then we have the following result.

\inspr{4.1} Theorem \rm 
The pair $(AB,\Delta_\#)$ is a regular multiplier Hopf algebra. The counit $\varepsilon_\#$ on $AB$ is given by
$$\varepsilon_\#(ab)=\varepsilon_A(a)\varepsilon_B(b)$$
and the antipode $S_\#$ is given by 
$$S_\#(ab)=\sum_{(a)}S_B(b)S_B(a_{(-1)})S_A(a_{(0)})$$
when $a\in A$ and $b\in B$.\hfill$\square$
\einspr

Remark that $a_{(-1)}$ is covered by $b$ in the last formula above.
\nl
It turns out that the easiest way to prove this result is by obtaining first the expressions for the linear maps $T_1^\#$ and $T_2^\#$, defined on $AB\ot AB$ as in the following proposition.

\inspr{4.2} Proposition \rm
Consider the linear maps $T_1^\#$ and $T_2^\#$ on $AB\ot AB$, defined by
$$\align T_1^\#(x\ot x')&=\Delta_\#(x)(1\ot x') \\
         T_2^\#(x\ot x')&=(x\ot 1)\Delta_\#(x').
\endalign$$
Then
$$\align T_1^\# &= (T^{-1})_{12}(T_1^A)_{23}T_{12}R_{34}(T_1^B)_{23}(R^{-1})_{34} \\
         T_2^\# &= (T^{-1})_{34}(T_2^B)_{23}T_{34}R_{12}(T_2^A)_{23}(R^{-1})_{12}.
\endalign$$
\hfill $\square$
\einspr

Before we give a proof of this important result, we need to make some remarks.
\snl
We use the leg numbering notation as explained before. The formulas not only involve the maps $T$ and $R$ but also the maps $T_1^A, T_2^A$ and $T_1^B, T_2^B$, defined like $T_1^\#, T_2^\#$, but for the multiplier Hopf algebras $A$ and $B$ respectively. 
\snl
Moreover, we must observe that the maps $T_1^\#$ and $T_2^\#$ are linear maps from  $AB\ot AB$ to itself, whereas on the right hand side of the equations, we have linear maps from $A\ot B\ot A\ot B$ to itself. Also, only the total expression will leave this space invariant. With the successive operations, we sometimes shuffle the tensor products. The reader may verify that indeed, the whole expression maps $A\ot B\ot A\ot B$ to itself, in the two cases. However, twists only occur within the first two factors or within the last two factors. For this reason, we can safely identify $A\ot B$ by $AB$ by means of the obvious map $a\ot b\mapsto ab$ within the first two and the last two factors. And because of this, we can use also that $ba=R(b\ot a)$. This observation is important for the proof below.
\snl
These remarks are sufficient for understanding and proving the result. We will give some more remarks of a different kind later, after the proof.

\inspr{} Proof \rm (of Proposition 4.2):
Take $a,a'\in A$ and $b,b'\in B$ and let $x=ab$ and $x'=b'a'$. The reader should be aware of the different order used to define $x$ and $x'$. 
\snl
By the remarks above, we have
$$(R_{34}(T_1^B)_{23}(R^{-1})_{34})(ab\ot b'a')=\sum_{(b)}ab_{(1)}\ot b_{(2)}b'a'.$$
Next, we want to apply the combination $(T_1^A)_{23}T_{12}$. This is a bit more tricky. First, consider $q,q'\in B$. Then
$$\align ((T_1^A)_{23}T_{12})(a\ot q\ot q'a')
             &=(T_1^A)_{23}(\Gamma(a)(q\ot 1)\ot q'a')\\
             &=((\iota_B\ot\Delta_A)\Gamma(a))(q\ot 1\ot q'a').
\endalign$$
where the right hand side is seen in $B\ot A\ot AB=B\ot A\ot BA$.
If we now replace $q\ot q'$ by $\Delta(b)(1\ot b')$, we arrive at
$$((T_1^A)_{23}T_{12}R_{34}(T_1^B)_{23}(R^{-1})_{34})(ab\ot b'a')=
               \sum_{(b)}((\iota_B\ot\Delta_A)\Gamma(a))(b_{(1)}\ot 1\ot b_{(2)}b'a').$$
\snl
Let us now look at $(T_{12} T_1^\#)(ab\ot b'a')$. We get 
$$\align (T_{12} T_1^\#)(ab\ot b'a')&=\sum_{(a)(b)}T_{12}((a_{(1)}\ot \Gamma(a_{(2)})(1\ot b_{(1)}\ot b_{(2)}b'a'))\\
               &= \sum_{(a)(b)}\Gamma_{12}(a_{(1)})\Gamma_{13}(a_{(2)})(b_{(1)}\ot 1\ot b_{(2)}b'a'),
\endalign$$ 
again seen in $B\ot A\ot AB$. 
\snl
This proves the first equality of the proposition as 
$$(\iota_B\ot \Delta_A)\Gamma(a)=\Gamma_{12}(a_{(1)})\Gamma_{13}(a_{(2)}).$$
\snl
The other equality is proven in a similar way, now taking $x=ba$ and $x'=a'b'$. Also, in the last step of the argument, now we have to use that
$$ (\Delta_B\ot\iota_A)\Gamma(a)=(\iota_B\ot\Gamma)\Gamma(a).$$
\hfill $\square$
\einspr

A reference to Majid's paper on the Hopf von Neumann algebra bicrossproducts ([M2]) is certainly appropriate here. It is well-known that the operators of the type $T_1$ are the algebraic counterparts of the fundamental unitaries $W$ as they appear in the theory of locally compact quantum groups (when working with the {\it right} Haar measure). The formula for $T_1^\#$, given in Proposition  4.2, is most easily recognized in the formula for $W$ as found in Exercise 6.2.14 of [M3]. The formula in Theorem 2.6 of [M2] looks a bit different (although still very similar). The reason for this difference is a consequence of the difference in conventions. In the theory of Kac algebra (or more generally, the locally compact quantum groups), the fundamental unitary $W$ is constructed from the left Haar measure. The formulas should also be compared with those in Definition 2.2 of [V-V2], but there the context is still more general and the comparison even more difficult.
\snl
As already mentioned in the introduction, these general analytical results do not imply our results because the setting is different. Nevertheless, it is instructive to compare the results.
\snl
In an earlier, unpublished version of this paper, our approach was different. The maps $T_1^\#$ and $T_2^\#$, as given in the proposition, were used to define the coproduct $\Delta_\#$. This has certain advantages. Coassociativity e.g.\ is proven by a straightforward verification of the equality
$$(T^\#_2\ot\iota)(\iota\ot T^\#_1)=(\iota\ot T^\#_1)(T^\#_2\ot\iota).$$

Also, that approach is more like in the analytical theories. However, for the reasons already explained in the introduction, in this paper, we have chosen another way, closer in spirit to the Hopf algebra case.
\nl
Now, we come to the proof of Theorem 4.1.

\inspr{} Proof\rm: 
We know already from the previous section that $\Delta_\#$ is a non-degenerate homomorphism from $AB$ to $M(AB\ot AB)$ and that it is coassociative. Then, it follows easily from the results in Proposition 4.2 that the pair $(AB,\Delta_\#)$ is a multiplier Hopf algebra. Indeed, the maps $T_1^\#$ and $T_2^\#$ are clearly bijective as compositions of bijective maps.
\snl
It is quite straightforward to obtain that $\varepsilon_\#$, as defined on $AB$ by $\varepsilon_\#(a)=\varepsilon_A(a)$ and $\varepsilon_\#(b)=\varepsilon_B(b)$ when $a\in A$ and $b\in B$, is the counit. Indeed, for $a\in A$ we have e.g.\ 
$$\align (\varepsilon_\#\ot \iota)(\Delta_\#(a))
                &=\sum_{(a)}(\varepsilon_\#\ot \iota)(a_{(1)}\ot 1)\Gamma(a_{(2)}) \\  
                &=\sum_{(a)}\varepsilon_A(a_{(1)})(\varepsilon_B \ot 1)\Gamma(a_{(2)}) \\  
                &=(\varepsilon_B\ot \iota)\Gamma(a)=a.
\endalign$$
To prove that also $(\iota\ot\varepsilon_\#)\Delta_\#(a)=a$, we use that $(\iota \ot \varepsilon_A)\Gamma(a)=\varepsilon_A(a)1$. Remark however that strictly speaking, only one equation must be proven. The other follows by the uniqueness of the counit. The formulas for the counit on the $B$-part are more or less trivial because the coproduct on this part coincides with the original coproduct.
\snl
Let us now look at the antipode. Again, we do not have to worry about $B$. Take $a\in A$ and first calculate $(S_\#\ot\iota)\Delta_\#(a)$. We get
$$\align \sum_{(a)}(S_\#\ot\iota)(a_{(1)}\ot 1)\Gamma(a_{(2)}) 
     &= \sum_{(a)}S_\#(a_{(1)}a_{(2)(-1)})\ot a_{(2)(0)}\\
     &= \sum_{(a)}S_B(a_{(2)(-1)})S_B(a_{(1)(-1)})S_A(a_{(1)(0)})\ot a_{(2)(0)}\\
     &= \sum_{(a)}S_B(a_{(1)(-1)}a_{(2)(-1)})S_A(a_{(1)(0)})\ot a_{(2)(0)}
\endalign$$
Now we use that $\Gamma_{12}(a_{(1)})\Gamma_{13}(a_{(2)})$ is equal to $(\iota_B\ot \Delta_A)\Gamma(a)$ and we find
$$(S_\#\ot\iota)\Delta_\#(a)=\sum_{(a)}S_B(a_{(-1)})S_A(a_{(0)(1)})\ot a_{(0)(2)}.$$
Finally, we multiply the two factors of this tensor product and we use that $m(S\ot\iota)\Delta(a')=\varepsilon(a')1$ for all $a'\in A$ and the fact that $(\iota_B\ot\varepsilon_A)\Gamma(a)=\varepsilon_A(a)1_B$ to get
$$m_\#(S_\#\ot\iota)\Delta_\#(a)=\sum_{(a)}S_B(a_{(-1)})\varepsilon_A(a_{(0)})=\varepsilon_A(a)1.$$
\snl
The last thing to observe is that the map $S_\#$, seen as a linear map from $A\ot B$ to itself, is given by the formula
$S_\#=R(S_B\ot S_A)T$. Therefore, it is a bijection of $AB$ and it follows that we have a regular multiplier Hopf algebra.
\hfill $\square$
\einspr

As we have seen in the previous section, if $A$ and $B$ are multiplier Hopf $^*$-algebras and if the action and coaction are compatible with the $^*$-structures, the algebra $AB$ is a $^*$-algebra and $\Delta_\#$ is a $^*$-homomorphism. Therefore, $(AB,\Delta_\#)$ is now a multiplier Hopf $^*$-algebra.
\nl
In the following theorem, we consider briefly the case of algebraic quantum groups. We claim that the bicrossproduct of algebraic quantum groups is again an algebraic quantum group. In a forthcoming paper, we will treat various aspects of this case in more detail and we will give more results and formulas (see [De-VD-W]).

\inspr{4.3} Theorem \rm
Let $A$ and $B$ be algebraic quantum groups. With the assumptions and notations of before (and the compatibility relations), also $AB$ will be an algebraic quantum group. In fact, if $\psi_A$ and $\psi_B$ are right integrals on $A$ and $B$ respectively, then a right integral $\psi_\#$ on $AB$ is given by $\psi_\#(ab)=\psi_A(a)\psi_B(b)$ when $a\in A$ and $b\in B$. If $A$ and $B$ are $^*$-algebraic quantum groups with positive right integrals  $\psi_A$ and $\psi_B$, then $AB$ is also a $^*$-algebraic quantum group and the right integral $\psi_\#$ given above is positive.

\snl\bf Proof\rm: 
We will first prove that $\psi_\#$ is right invariant. 
\snl
Start with two elements $a,a''\in A$. We can write
$$\sum_{(a)}\psi(a_{(1)})\Delta_\#(a''a_{(2)})=\psi(a)\Delta_\#(a'').$$
We now let $\Delta_\#(A)$ act from the right on $B\ot A$ as before (see Notation 3.4). Then, for all $a,a',a''\in A$ and $b'\in B$, we have
$$\sum_{(a)}\psi(a_{(1)})((b'\ot a')\bullet\Delta_\#(a''a_{(2)}))=\psi(a)(b'\ot a')\bullet\Delta_\#(a'').$$
Because this action of $A$ is unital, it follows that
$$\sum_{(a)}\psi(a_{(1)})(b'\ot a')\bullet\Delta_\#(a_{(2)})=\psi(a)(b'\ot a'),$$
still for all $a,a'\in A$ and $b'\in B$. Observe that in this expression, $a_{(2)}$ is covered by $b'\ot a'$ through the action.
\snl
Now, take $a,a'\in A$ and $b,b'\in B$. Then we have
$$\align (\psi_\#\ot\iota)((1\ot a'b')\Delta_\#(ba))
         &=\sum_{(b)}(\psi_\#\ot\iota)((b_{(1)}\ot a'b'b_{(2)})\Delta_\#(a))\\
         &=\sum_{(a)(b)}(\psi_\#\ot\iota)((a_{(1)}\ot 1)((b_{(1)}\ot a'b'b_{(2)})\bullet\Delta_\#(a_{(2)})))\\
         &=\sum_{(a)(b)}\psi_A((a_{(1)})(\psi_B\ot\iota)((b_{(1)}\ot a'b'b_{(2)})\bullet\Delta_\#(a_{(2)})))\\
         &=\sum_{(b)}\psi_A(a)(\psi_B\ot\iota)((b_{(1)}\ot a'b'b_{(2)}))\\
         &=\psi_A(a)\psi_B(b)a'b'.
\endalign$$
Now, clearly $\psi_\#(ba)=\sum_{(a)}\psi_A(a_{(1)})\psi_B(b\tl a_{(2)})=\psi_A(a)\psi_B(b)$. This proves that $\psi_\#$ is indeed the right invariant integral on the smash product $AB$.
\snl
Now assume that $A$ and $B$ are $^*$-algebras and that $\psi_A$ and $\psi_B$ are positive. Then, for all $a,a'\in A$ and $b,b'\in B$ we have
$$\align \psi_\#((ab)^*(a'b'))&=\psi_\#(b^*a^*a'b') \\
             &=\sum_{(a^*)(a')}\psi_A(a^*_{(1)}a'_{(1)})\psi_B((b^*\tl(a^*_{(2)}a'_{(2)}))b')\\
             &=\psi_A(a^*a')\psi_B(b^*b').
\endalign$$
It follows that also $\psi_\#$ will be positive. \hfill $\square$
\einspr

When $A$ has an identity, the argument is much easier because one can use that $\Gamma(1)=1$ in $M(B\ot A)$. In some sense, this is still true as, by general results on multiplier Hopf algebras, we can extend the coproduct $\Delta_\#$ to the multiplier algebra and we get $\Delta_\#(1)=1$. This is essentially the same as $\Gamma(1)=1$. However, because $\Gamma$ is not a homomorphism, we can not simply apply the rules to extend it to the multiplier algebra.
\snl
In our forthcoming paper on algebraic quantum groups [De-VD-W], we will also obtain expressions for the left integral on the bicrossproduct and the other data associated with an algebraic quantum group and its dual. Related results are obtained in [B-*] and in [De-VD-W] we will discuss further the relation between these different approaches.
\nl
Next, we consider the dual case. So, we have two regular multiplier Hopf algebras $C$ and $D$, but now we assume that $C$ is a left $D$-module algebra and that $D$ is a right $C$-comodule coalgebra. We use the notations as recalled in the previous section and we assume the conditions as formulated in Theorem 3.16. Then we get the following analogue of Theorem 4.1.

\inspr{4.4} Theorem \rm 
The pair $(CD,\Delta_\#)$ is a regular multiplier Hopf algebra. The counit $\varepsilon_\#$ on $CD$ is given by
$$\varepsilon_\#(cd)=\varepsilon_C(c)\varepsilon_D(d)$$
and the antipode $S_\#$ is given by
$$S_\#(cd)=\sum_{(d)} S_D(d_{(0)})S_C(d_{(1)})S_C(c)$$
whenever $c\in C$ and $d\in D$. \hfill$\square$

\einspr

There is no need to prove this result as it follows from Theorem 4.1 and the technique to pass from the pair $(A,B)$ to the pair $(C,D)$ as explained earlier in this paper. When $A$ is paired with $C$ and when $B$ is paired with $D$ and if the actions and coactions are dual to each other, we get a pairing of the regular multiplier Hopf algebras $AB$ and $CD$ as explained in the previous section. In the $^*$-algebra case, we get a pairing of multiplier Hopf $^*$-algebras. 
\snl
We also have the analogue of Theorem 4.3. If the algebras $C$ and $D$ are algebraic quantum groups, then $CD$ is again an algebraic quantum group and a left integral $\varphi_\#$ on $CD$ is given by $\varphi_\#(cd)=\varphi_C(c)\varphi_D(d)$ where $\varphi_C$ and $\varphi_D$ are left integrals on $C$ and $D$ respectively. Also, if the algebras are $^*$-algebraic quantum groups with positive integrals, then this is the case for $CD$.
\snl
In our paper [De-VD-W], we will show that, in the case of algebraic quantum groups, and when $C$ and $D$ are the duals $\hat A$ of $A$ and $\hat B$ of $B$ respectively, the smash product $CD$ is naturally identified with the dual $(AB)^{\hat{}}$ of the smash product $AB$. The result is not completely obvious as will be seen in that paper.
\nl
There is nothing more to say about the special cases. The important results with references to the literature have been treated in the previous section. Similarly, not much more can be said about the examples. We just look at the formulas for the counit and the antipode, found in Theorems 4.1 and 4.4.

\inspr{4.5} Examples \rm
Consider the examples in 3.17. In the first case, the counit $\varepsilon_\#$ on $AB$ is given by $\varepsilon_\#(hf)=\varepsilon_A(h)\varepsilon_B(f)=f(e)$ for $h\in H$ and $f\in F(K)$ where $e$ is the unit of the group. Similarly, the counit on $CD$ is given by the formula $\varepsilon_\#(fk)=f(e)$ for $k\in K$ and $f\in F(H)$. 
\snl
For the antipode on $AB$ we get
$$S_\#(h\delta_k)=S_B(\delta_k)S_A(h\tl k)=\delta_{k^{-1}}(h\tl k)^{-1}$$
when $h\in H$ and $k\in K$. A careful calculation gives that this last expression is also equal to $(h\tl k)^{-1}\delta_{(h\tr k)^{-1}}$. Compare with the formula in Example 6.2.12 of [M3]]. Similarly, we find for the antipode on $CD$:
$$S_\#(\delta_h k)= \delta_{(h\tl k)^{-1}}(h\tr k)^{-1}$$
and this should be compared with the formula in Example 6.2.11 of [M3].
\hfill $\square$
\einspr

\nl\nl

\bf 5. Conclusion and further research \rm
\nl
We have shown in this paper how the bicrossproduct construction of Majid can be generalized to regular multiplier Hopf algebras. The main result is Theorem 4.1 in Section 4.
\snl
We started with two regular multiplier Hopf algebras $A$ and $B$, a right action $\tl$ of $A$ on $B$ making $B$ into a right $A$-module algebra (see Section 1), a left coaction $\Gamma$ of $B$ on $A$ making $A$ into a left $B$-comodule coalgebra (see Section 2) and natural compatibility relations between the action and the coaction (see Section 3). We have also treated the dual case.
\snl
The notion of a right module algebra and the associated smash product, reviewed in Section 1, presents no problem and is well understood. It has been studied for multiplier Hopf ($^*$-)algebras in [Dr-VD-Z]. The aspects of coverings are relatively simple ones and are used in [VD4] to illustrate a new (and better) approach to the technique of covering.
\snl
The situation is somewhat more complicated when it comes to coactions and smash coproducts. This has been discussed in Section 2. The treatment is slightly different and in particular, more general (and in some sense also more natural) than the one in earlier literature (see e.g.\ [VD-Z2] and [De2]). Many aspects of covering are more complicated. To give a meaning to the basic requirement
$$(\iota_B \ot \Delta_A)\Gamma(a)=\sum_{(a)}\Gamma_{12}(a_{(1)})\Gamma_{13}(a_{(2)})$$
for a left $B$-comodule coalgebra $A$ with coaction $\Gamma$ and of course $a\in A$, one needs an extra condition (see Assumption 2.6). However, as mentioned in Section 2, this does not completely fit into the covering framework as treated in the appendices. This is unfortunate.
\snl 
Therefore, it seems that some more research is needed. A natural condition to impose on $\Gamma$ would be e.g.\ that the variable $p$ of $A$ in the expression $\Gamma(p)(b\ot a)$ is covered when $a\in A$ and $b\in B$ are given. This, together with similar conditions,  would imply the assumption needed in Section 2. There are also reasons to believe that this is indeed the case. However, it is probably necessary to have the compatibility relations of Section 3, relating the coaction $\Gamma$ with the product on $A$ for the simple reason that covering a variable in $A$ only makes sense when $A$ is an algebra.
\snl
Now, we come to Section 3 where we study the necessary conditions that relate the action with the coaction. First, there is the condition needed for the coproduct $\Delta_\#$ to be multiplicative on $A$. It shows how to express $\Gamma(aa')$ when $a,a'\in A$ in terms of $\Gamma(a)$ and $\Gamma(a')$. We have this condition under a nice form in terms of the map $P$ (see Proposition 3.3). Next, we have the natural dual version of this condition, as formulated in Proposition 3.10. It relates the coproduct of $B$ with the action of $A$ on $B$ and expresses $\Delta_B(b\tl a)$ in terms of $\Delta_B(b)$ and $\Gamma(a)$. Again, we have a form of this condition in terms of the map $P$ (cf.\ Assumption 3.9). Finally, there is the (somewhat strange) equality $T\circ R=T^{\text{op}}\circ R^{\text{op}}$ (Assumption 3.12), needed together with the other assumptions, to have that the natural coproduct $\Delta_\#$ is indeed an algebra map on the smash product $AB$.
\snl
There are no real extra difficulties when it comes to covering the formulas in Section 3. There is only the condition 
$$(\Delta_B\ot \iota_A)P=P_{23}P_{13}(\Delta_B\ot \iota_A)$$
that needs an extra, but simple argument.
\snl
The conditions, as formulated by us in this Section 3, are of course essentially the same as the ones you get when generalizing the original conditions of Majid to the case of multiplier Hopf algebras. We wonder however if they have been considered earlier in these other forms for Hopf algebras?
\snl
In Section 4, we have obtained the main result (Theorem 4.1). However, the 'hard work' was done already in Section 3. The difficulty that remains is to show in Proposition 4.2 that the maps $T_1^\#$ and $T_2^\#$, canonically associated with the coproduct $\Delta_\#$, can be given in terms of the other given linear maps. This then allows to show easily that these canonical maps are bijective and so that the pair $(AB,\Delta_\#)$ is a multiplier Hopf algebra. Regularity of this multiplier Hopf algebra is obtained by showing that the antipode $S_\#$ is bijective. From the regularity, it should be possible to answer some of the questions, raised earlier about covering properties in Section 2, with solutions depending upon the conditions of Section 3. Also here, some more investigations would be welcome. See also an earlier remark above.  
\snl
We have given a  rigorous proof of Theorem 4.3 where it is shown that the bicrossproduct of algebraic quantum groups is again an algebraic quantum group and where a formula for the right integral on $AB$ is given. The argument seems more complicated than when $A$ has an identity. What is really needed  is an interpretation and an argument for the formula $\Gamma(1_A)=1_B\ot 1_A$. The problem is of course related with the problems formulated before about aspects of covering for the coaction. So, again, more research would help for better understanding these relations. We intend to consider this in our forthcoming paper on bicrossproducts of algebraic quantum groups ([(D-VD-W]).  
\nl
In this second paper on the subject, we do not only treat the integrals on the bicrossproduct of $A$ and $B$ (when $A$ and $B$ have integrals). We also show that the duality, considered in this paper between $AB$ and $CD$, realizes $CD$ as the dual of $AB$ (in the sense of algebraic quantum groups), when $C=\hat A$ and $D=\hat B$ and when of course the actions and coactions are adjoints of each other. This allows to formulate many relations between the various objects associated with integrals on multiplier Hopf algebras (such as the modular element, the modular automorphism groups and the scaling constants).
\nl
Finally, there is the need to construct examples, other than the one given in this paper coming from the decomposition of a group $G$ into two subgroups $H$ and $K$. Of course, we want examples with multiplier Hopf algebras, if possible with integrals, that are not ordinary Hopf algebras. There have been constructed new and interesting examples of the bicrossproduct construction for the more general locally compact quantum groups (see [V1], [V2], [V-V1] and [V-V2]), but there is little hope that these will fit into this framework here.
\snl
On the other hand, it is now well understood when a locally compact group $G$ gives rise to a multiplier Hopf algebra. This is the case if and only if the group has a compact open subgroup (see [L-VD]). Therefore, it is expected that examples of bicrossproducts in our framework of multiplier Hopf algebras (with integrals), will be found using such groups. It is also conceivable that some of the examples of Vaes and Vainerman can be modified so as to fit into our treatment. We hope to do this also in our next paper.
\nl\nl

\bf References \rm
\nl
{\bf [B-S]} S.\ Baaj and G.\ Skandalis: {\it Unitaires multiplicatifs et dualit\'e pour les produits crois\'es de C$^*$-alg\`ebres}. Ann.\ Scient.\ Ec.\ Norm.\ Sup.\ 4$^{\text{e}}$ s\'erie {\bf 26} (1993), 425--488. 
\snl
{\bf [B-*]} M.\ Beattie, S.\ D\v{a}sc\v{a}lescu, L.\ Gr\"unenfelder and C.\ N\v{a}st\v{a}sescu: {\it Finiteness conditions, co-Frobenius Hopf algebras and quantum groups}. J.\ Algebra {\bf 200} (1998), 312-333.
\snl
{\bf [B-M]} T.\ Brzezinski and S.\ Majid: {\it Quantum Geometry of Algebra Factorisations and Coalgebra Bundles}. Commun.\ Math.\ Phys.\ {\bf 213} (2000), 491-521.
\snl
{\bf [De1]} L.\ Delvaux: {\it Semi-direct products of multiplier Hopf algebras: Smash products}. Comm. Algebra {\bf 30} (2002), 5961--5977.
\snl
{\bf [De2]} L.\ Delvaux: {\it Semi-direct products of multiplier Hopf algebras: Smash coproducts}. Comm. Algebra {\bf 30} (2002), 5979-5997.
\snl
{\bf [De3]} L.\ Delvaux: {\it Twisted tensor product of multiplier Hopf ($^\ast$-)algebras}. J.\ Algebra {\bf 269} (2003), 285-316.
\snl
{\bf [De4]} L.\ Delvaux: {\it Twisted tensor coproduct of multiplier Hopf ($^\ast$-)algebras}. J.\ Algebra {\bf 274} (2004), 751--771.
\snl
{\bf [Dr-VD]} B.\ Drabant \& A. Van Daele: {\it Pairing and Quantum double of multiplier Hopf algebras}.  Algebras and Representation Theory 4 (2001), 109-132.
\snl
{\bf [De-VD-W]} L.\ Delvaux, A.\ Van Daele \& S.\ Wang: {\it Bicrossproducts of algebraic quantum groups}. Preprint U.\ Hasselt, K.U.\ Leuven \& Nanjing University. In preparation.
\snl
{\bf [Dr-VD-Z]} B.\ Drabant, A.\ Van Daele \& Y.\ Zhang: {\it Actions of multiplier Hopf algebras}. Comm.\ Algebra {\bf 27} (1999), 4117-4127.
\snl
{\bf [L-VD]} M.B.\ Landstad \& A.\ Van Daele: {\it Groups with compact open subgroups and multiplier Hopf $^*$-algebras}. Expo.\ Math.\ 26 (2008), 197--217. 
\snl
{\bf [M1]} S.\ Majid: {\it Physics for algebraists: Non-commutative and non-cocommutative Hopf algebras by a bicrossproduct construction}. J.\ Algebra {\bf 130} (1990), 17-64.
\snl
{\bf [M2]} S.\ Majid:  {\it Hopf-von Neumann algebra bicrossproducts, Kac algebra bicrossproducts and the classical Yang-Baxter equations}. J.\ Funct.\ Analysis {\bf 95} (1991), 291-319.
\snl
{\bf [M3]} S.\ Majid: {\it Foundations of quantum group theory}. Cambridge University Press, 1995.
\snl
{\bf [V1]} S.\ Vaes: {\it Locally compact quantum groups}. Ph.D. thesis K.U.Leuven (2001).
\snl
{\bf [V2]} S.\ Vaes: {\it Examples of locally compact quantum groups through the bicrossed product construction}. Proceedings of the XIIIth Int.\ Conf.\ Math.\ Phys.\ London (2000). Editors A.\ Grigoryan, A.\ Fokas, T.\ Kibble and B.\ Zegarlinski, International press of Boston, Somerville MA (2001), pp.341--348.
\snl
{\bf [V-V1]} S.\ Vaes and L.\ Vainerman: {\it On low dimensional locally compact quantum groups}. Locally Compact Quantum Groups and Groupoids. Proceedings of the Meeting of Theoretical Physicists and Mathematicians, Strasbourg (2002). Ed.\ L.\ Vainerman, IRMA Lectures on Mathematics and Mathematical Physics, Walter de Gruyter, Berlin, New York (2003), pp.127-187.
\snl
{\bf [V-V2]} S.\ Vaes and L.\ Vainerman: {\it Extensions of locally compact quantum groups and the bicrossed product construction}. Adv.\ in Math.\ {\bf 175} (2003), 1--101.  
\snl
{\bf [VD1]} A.\ Van Daele: {\it Multiplier Hopf algebras}. Trans.\ Amer.\ Math.\ Soc.\ {\bf 342} (1994), 917-932.
\snl
{\bf [VD2]} A.\ Van Daele: {\it An algebraic framework for group duality}. Adv.\ Math.\ {\bf 140} (1998), 323-366.
\snl
{\bf [VD3]} A.\ Van Daele: {\it Multiplier Hopf $^*$-algebras with positive integrals: A laboratory for locally compact quantum groups}. Locally Compact Quantum Groups and Groupoids. Proceedings of the Meeting of Theoretical Physicists and Mathematicians, Strasbourg (2002). Ed.\ L.\ Vainerman, IRMA Lectures on Mathematics and Mathematical Physics 2, Walter de Gruyter, Berlin, New York (2003), pp.229-247.
\snl
{\bf [VD4]} A.\ Van Daele: {\it Tools for working with multiplier Hopf algebras}. Arabian Journal for Science and Engineering {\bf 33} 2C (2008), 505--527. See also Arxiv math.RA/0806.2089. 
\snl
{\bf [VD-VK]} A.\ Van Daele \& S.\ Van Keer: {\it The Yang-Baxter and the Pentagon equation}. Comp.\ Math.\ {\bf 91} (1994), 201--221.
\snl 
{\bf [VD-W]} A.\ Van Daele \& S.\ Wang: {\it A class of multiplier Hopf algebras}. Algebras and Representation Theory {\bf 10} (2007), 441--461.
\snl
{\bf [VD-Z1]} A.\ Van Daele \& Y.\ Zhang: {\it A survey on multiplier Hopf algebras} In 'Hopf algebras and Quantum Groups', eds. S.\ Caenepeel \& F.\ Van Oyestayen, Dekker, New York (1998), pp. 259--309.
\snl
{\bf [VD-Z2]} A.\ Van Daele and Y.\ Zhang: {\it Galois Theory for multiplier Hopf algebras with integrals}. Algebras and Representation Theory {\bf 2} (1999), 83-106.

\end